\documentclass[10pt]{article}

\usepackage{Variationen}

\author{II. Symmetric powers of $\Nat \sl_2(\K)$}

\begin{document}

\setcounter{variation}{13}

\renewcommand{\fg}{\mathfrak{g}}

\maketitle

\liminaryquotation{\flushright C'est une série logarithmique.}

\abstract{We identify the spaces of homogeneous polynomials in two variables $\K[Y^k, XY^{k-1}, \dots, X^k]$ among representations of the Lie ring $\sl_2(\K)$. This amounts to constructing a compatible $\K$-linear structure on some abstract $\sl_2(\K)$-modules, where $\sl_2(\K)$ is viewed as a Lie ring.}

\renewcommand{\thefootnote}{}

The present article comes immediately next to \cite{TV-I} but may be read independently. One characterizes the symmetric powers $\Sym^k \Nat \sl_2(\K)$ among representations of the Lie ring $\sl_2(\K)$.

Remember that the modules $\Sym^k \Nat \sl_2(\K)$ ($k \geq 1$ an integer) are the irreducible, finite-dimensional representations of the Lie $\K$-algebra $\sl_2(\K)$. It is convenient to realize them as the spaces of homogeneous polynomials of degree $k$ in two variables $\K[Y^k, XY^{k-1}, \dots, X^k]$ equipped with the standard action.

Such are the irreducible finite-dimensional representations of the Lie $\K$-\emph{algebra}. But when one sees $\sl_2(\K)$ as a mere Lie \emph{ring}, the question is more general. For one simply lets $\sl_2(\K)$ act on an arbitrary abelian group $V$ in such a way that $\sl_2(\K) \to \End(V)$ be a morphism of Lie rings. In this abstract setting $V$ need not be a $\K$-vector space, which leaves us quite far from weight theory. The present article thus deals like \cite{TV-I} with linearizing abstract modules, that is constructing a $\K$-vector space structure compatible with the given action of an algebraic structure, here the Lie ring $\sl_2(\K)$.

It is not surprising to turn the assumption that $V$ is finite-dimensional (which is a priori meaningless since there is no $\K$-linear structure to start with) into an assumption on the length of the action: $x^n \cdot V = 0$, where $h, x, y$ form the standard $\K$-basis of $\sl_2(\K)$. We thus extend the results of \cite{TV-I} which considered the simpler, quadratic case $x^2 \cdot V = 0$. Actually much more information on $V$ is given by its description as an $\sl_2(\K_1)$-module where $\K_1$ is the prime subfield of $\K$. Whenever we write $A \simeq \oplus_I B$ we merely mean that $A$ is isomorphic to a direct sum of copies of $B$ indexed by some (possibly finite) set $I$. Our main result is the following.\medskip

\par\noindent\textbf{Variation \ref{v:Symn-1Natsl2K}.}
{\itshape
Let $n \geq 2$ be an integer and $\K$ be a field of characteristic $0$ or $\geq n$.
Let $\fg = \sl_2(\K)$ viewed as a Lie ring and $V$ be a $\fg$-module. Let $\K_1$ be the prime subfield of $\K$ and $\fg_1 = \sl_2(\K_1)$. Suppose that $V$ is a $\K_1$-vector space such that $V \simeq \oplus_I \Sym^{n-1} \Nat \fg_1$ as $\K_1\fg_1$-modules.

Then $V$ bears a compatible $\K$-vector space structure for which $V \simeq \oplus_J \Sym^{n-1} \Nat \fg$ as $\K\fg$-modules.}
\medskip

One may actually say a little more under an extra hypothesis which we shall call of \emph{coherence} of the action, in the sense that the kernels and/or images of the nilpotent operators must obey some global behaviour. 
\medskip

\par\noindent\textbf{Variations \ref{v:suiteSymNatsl2K:ker}, \ref{v:suiteSymNatsl2K:im} and \ref{v:oplusSymNatsl2K}.}
\textit{
Let $n \geq 2$ be an integer and $\K$ be a field.
Let $\fg = \sl_2(\K)$ viewed as a Lie ring and $V$ be a $\fg$-module. If the characteristic of $\K$ is $0$ one requires $V$ to be torsion-free.
Suppose that:
\begin{itemize}
\item
either $x^n = 0$ in $\End V$ and the characteristic of $\K$ is $0$ or $\geq 2n+1$,
\item
or $x^n = y^n = 0$ in $\End V$ and the characteristic of $\K$ is $\geq n+1$.
\end{itemize}
Then:
\begin{itemize}
\item
if for all $\lambda$, $\ker x \leq \ker x_\lambda$, then there is a series $\Ann_V(\fg) = V_0 \leq V_1 \leq \dots \leq V_{n-1} = V$ of $\fg$-submodules such that for all $k = 1, \dots, n-1$, $V_k/V_{k-1}$ bears a compatible $\K$-vector space structure for which $V_k/V_{k-1} \simeq \oplus_{I_k} \Sym^k \Nat \fg$ as $\K\fg$-modules (\ref{v:suiteSymNatsl2K:ker});
\item
if for all $\lambda \in \K$, $\im x_\lambda \leq \im x$, then there is a series $0 = V_0 \leq V_1 \leq \dots \leq V_{n-1} = \fg \cdot V$ of $\fg$-submodules such that for all $k = 1, \dots, n-1$, $V_k/V_{k-1}$ bears a compatible $\K$-vector space structure for which $V_k/V_{k-1} \simeq \oplus_{I_{n-k}} \Sym^{n-k} \Nat \fg$ as $\K\fg$-modules (\ref{v:suiteSymNatsl2K:im});
\item
if for all $\lambda$, $\ker x \leq \ker x_\lambda$ \emph{and} $\im x_\lambda \leq \im x$, then our series split: $V = \Ann_V(\fg)\oplus \fg \cdot V $, and $\fg\cdot V$ bears a compatible $\K$-vector space structure for which $\fg \cdot V \simeq \oplus_{k = 1}^{n-1} \oplus_{I_k} \Sym^k \Nat \fg$ as $\K\fg$-modules (\ref{v:oplusSymNatsl2K}).
\end{itemize}
}
\medskip

One should in particular note the immediate consequence:
\begin{quote}
if $V$ is a \emph{simple} $\fg$-module with $n$ minimal such that $x^n \cdot V = 0$, and -- either for all $\lambda \in \K$, $\ker x \leq \ker x_\lambda$ -- or for all $\lambda \in \K$, $\im x_\lambda \leq \im x$, then $V \simeq \Sym^{n-1} \Nat \fg$.
\end{quote}
Without simplicity the statement we gave is quite clumsier due to the trouble one has controlling the cohomology of representations of a Lie ring. The reader will observe that we cannot prove in general that $\fg\cdot V/\Ann_{\fg\cdot V} \fg$ bears a compatible $\K$-vector space structure under either assumption $\ker x \leq \ker x_\lambda$ or $\im x_\lambda \leq \im x$: we seem to need both.

The present article starts with a few notations (\S\ref{S:cadre}) and basic remarks on length (\S\ref{S:length}). Section \ref{S:squelette} then studies the actions of $\sl_2(\K_1)$ for a \emph{prime} field $\K_1$. Everything is as expected in large enough characteristic (\S\ref{s:pgeq2n+1}).
In \S\ref{s:pathologies} we somehow digress by lowering a little the characteristic which results in creating pathologies. These are removed in \S\ref{s:xn=yn=0} under the assumption that the action is decently ``two-sided''. At this point we leave prime fields for the general case and move to Section \ref{S:chair}.
Our main result Variation \ref{v:Symn-1Natsl2K} describes the extension of the linear structure from the combinatorial skeleton (i.e., the action at the level of the prime subfield $\K_1$) to the scalar flesh (i.e., the action at the level of the field $\K$); it is proved in \S\ref{s:separe}. In \S\ref{s:suites} the inclusions $\ker x \leq \ker x_\lambda$ and $\im x_\lambda \leq \im x$ finally appear. A few ideas on their possible meaning are put forth in Section \ref{S:coherence}.

Technically speaking the only tool is patience; the reader should expect long computations. Studying $\sl_2(\K)$ as a Lie ring may sound somehow arbitrary, and the author has no illusions on his results. The main purpose of the study was to prepare him for the future case of the group $\SL_2(\K)$.

This work was finished during a visit to the French-Russian ``Poncelet'' Mathematics Laboratory in Moscow. The author warmly thanks everyone involved, with a special thought for the gentleman who likes Belgian chocolates.

\section{The setting}\label{S:cadre}

This section contains notations and very basic facts which will be used with no reference.

\subsection{The Lie ring}

\begin{notation*}
Let $\K$\inmargin{$\K$, $\fg$} be a field and $\fg$ be the Lie ring $\sl_2(\K)$.
\end{notation*}

Literature on Lie rings looks scarce when compared to other topics. Fortunately we deal with a concrete and familiar Lie ring, so any reference on Lie algebras such as \cite{Humphreys} will do. We simply use the group law $+$, the bracket $[\cdot, +\cdot]$, and forget about the $\K$-linear structure on $\fg$.

\begin{notation*}
Let $\K_1$\inmargin{$\K_1$, $\fg_1$} be the prime subfield of $\K$ and $\fg_1$ be the Lie ring $\sl_2(\K_1)$; one has $\fg_1 \leq \fg$.
\end{notation*}

\begin{notation*}
For $\lambda \in \K$ let\inmargin{$h_\lambda$, $x_\lambda$, $y_\lambda$}
\[h_\lambda = \left(
\begin{array}{cc}
 \lambda & 0 \\ 0 & -\lambda
\end{array}\right), \quad
x_\lambda = \left(
\begin{array}{cc}
 0 & \lambda \\ 0 & 0
\end{array}\right), \quad
y_\lambda = \left(
\begin{array}{cc}
 0 & 0 \\ \lambda & 0
\end{array}\right)\]
One simply writes $h = h_1$, $x = x_1$, $y = y_1$.\inmargin{$h, x, y$}
\end{notation*}

\begin{notation*}
Let $\fb$\inmargin{$\fb, \fu, \ft$} be the Borel subring generated by the $h_\lambda$'s and the $x_\mu$'s, and $\fu = \{x_\mu: \mu \in \K\}$ be its unipotent radical. Let $\ft$ be the Cartan subring $\{h_\lambda : \lambda \in \K\}$.
\end{notation*}

\begin{relations*}\
\begin{itemize}
\item
$[h_\lambda, x_\mu] = 2 x_{\lambda \mu}$;
\item
$[h_\lambda, y_\nu] = - 2 y_{\lambda \nu}$;
\item
$[x_\mu, y_\nu] = h_{\mu\nu}$.
\end{itemize}
\end{relations*}

$\K$ will never have characteristic $2$; as a consequence $\fg$ will always be perfect. One should be careful that $[\fg, \fg]$ will merely denote the \emph{additive subgroup} of $\fg$ generated by all brackets since we forget about the $\K$-linear structure on $\fg$. It is however the case that $\fg = [\fg, \fg]$ which is the definition of perfectness.

We shall sometimes go to the enveloping (associative) ring which is defined among rings just like the enveloping (associative) $\K$-algebra is defined among $\K$-algebras. It enjoys a similar universal property in the broader category of $\fg$-modules. This simply amounts to viewing $\sl_2(\K)$ as a Lie algebra over the prime ring of $\K$ and taking its enveloping algebra as such, but to prevent confusion we shall always refer to this object as the enveloping ring. The usual enveloping $\K$-algebra can be retrieved as a quotient of the enveloping Lie ring by relations expressing $\K$-linearity. It has no reason to play a role since we are a priori not dealing with $\K$-linear representations.

\begin{relations*}
One has in the enveloping ring the following equalities which the reader may check by induction:
\begin{align}
x^i y & = y x^i + i (h+1-i) x^{i-1};\\
y^j x & = xy^j - j (h+j-1)y^{j-1};\\
y^jx^i & = \sum_{k = 0}^{\min(i, j)} \left[(-1)^k k!\ \binom{i}{k} \binom{j}{k} \left(\prod_{\ell = 0}^{k-1}(h-i+j+\ell)\right) x^{i-k}y^{j-k} \right];\\
x^i h_\mu & = h_\mu x^i - 2i x_\mu x^{i-1};\\
x_{\lambda_1} \dots x_{\lambda_i} y_\mu & = \quad y_\mu x_{\lambda_1} \dots x_{\lambda_i} + \sum_{k } h_{\mu \cdot \lambda_k} x_{\lambda_1} \dots \widehat{x_{\lambda_k}} \dots x_{\lambda_i}\nonumber\\
& \quad - \sum_{k \neq \ell} x_{\mu\cdot \lambda_k \cdot \lambda_\ell}x_{\lambda_1} \dots \widehat{x_{\lambda_k}} \dots \widehat{x_{\lambda_\ell}} \dots x_{\lambda_i};\\
y_\mu x_{\lambda_1} \dots x_{\lambda_i} & = \quad x_{\lambda_1} \dots x_{\lambda_i} y_\mu - \sum_{k } x_{\lambda_1} \dots \widehat{x_{\lambda_k}} \dots x_{\lambda_i} h_{\mu \cdot \lambda_k}\nonumber \\
& \quad - \sum_{k \neq \ell} x_{\lambda_1} \dots \widehat{x_{\lambda_k}} \dots \widehat{x_{\lambda_\ell}} \dots x_{\lambda_i} x_{\mu\cdot \lambda_k \cdot \lambda_\ell}
\end{align}
(The terms in the hats do not appear.)
\end{relations*}

Be however careful that in the enveloping ring $x_\lambda y \neq x y_\lambda$. So checking the formulas in the enveloping algebra does not suffice in order to establish them in the enveloping ring.

\begin{notation*}
Let $c_1 = 2 xy + 2yx + h^2$ be the Casimir operator.
\end{notation*}

The Casimir operator is central in the enveloping algebra but not in the enveloping ring; for instance a quick computation yields $[c_1, h_\lambda] = 8 (xy_\lambda - x_\lambda y)$ which is non-zero.

However when $\K_1 = \F_p$ and $\fg_1 = \sl_2(\K_1)$, $c_1$ is central indeed in the enveloping ring of $\fg_1$. This is not quite true over $\Q$, but it is readily checked that for all $z$ in the enveloping ring there is an integer $k$ with $k [c_1, z] = 0$. It follows that provided $\K_1 = \Q$ and $V$ is a \emph{torsion-free} $\fg_1 = \sl_2(\K_1)$-module, the action of $c_1$ commutes with the action of $\fg_1$. This will always be the case when we use $c_1$.

\subsection{The module}

\begin{notation*}
Let $V$\inmargin{$V$} be a $\fg$-module.
\end{notation*}

We shall keep writing $x_\lambda, y_\lambda, h_\mu$ for the images in $\End V$ of the corresponding elements of $\fg$.

\begin{notation*}
The length of the $\fu$-module\inmargin{$\lambda_\fu(V)$} $V$ is the least integer $n$, if there is one such, with $\fu^n \cdot V = 0$.
\end{notation*}

\begin{notation*}
For $i \in \Z$, let $E_i = E_i(V) = \{a \in V: h\cdot a = iv\}$\inmargin{$E_i$}.
\end{notation*}

Using the familiar relations one sees that $h_\lambda$ (resp. $x_\mu$, resp. $y_\nu$) maps $E_i$ to $E_i$ (resp. $E_{i+2}$, resp. $E_{i-2}$). 

We shall in a minute deal with constructing vector space structures on modules. If $\K_1$ is a prime field then an abelian group $V$ bears at most one structure over $\K_1$. If it is the case and $V$ is a $\fg_1 = \sl_2(\K_1)$-module as well then $V$ is a $\K_1\fg_1$-module.

Let us also remind the reader why Lie rings do not admit cross-characteristic representations.

\begin{observation*}
Let $\K$ be a field, $\K_1$ its prime subfield, $\fg = \sl_2(\K)$ and $V$ be a $\fg$-module. Then $\fg\cdot V/\Ann_{\fg\cdot V}\fg$ is a $\K_1$-vector space.
\end{observation*}
\begin{proofclaim}
If $\K_1 = \F_p$ then $\fg$ annihilates $pV$ so $V/\Ann_V \fg$ has exponent $p$. Also note that $p$ annihilates $\fg\cdot V$, so $\fg\cdot V$ has exponent $p$ as well. Hence in prime characteristic both $V/\Ann_V \fg$ and $\fg \cdot V$ are actually $\K_1$-vector spaces.

If $\K_1 = \Q$ then $\fg\cdot V$ is divisible and $\fg$ annihilates the torsion submodule of $V$, so $\fg\cdot V/\Ann_{\fg\cdot V}\fg$ is torsion-free and divisible: a $\Q$-vector space.
\end{proofclaim}

This certainly does not prove that $V$ need be a $\K$-vector space (which is not true in general) but already removes the outmost pathologies.

\subsection{Symmetric powers}

Let us finally recast a few facts about the very modules we try to characterize.

\begin{notation*}
Let $\Nat \fg$\inmargin{$\Nat \fg$} denote the natural representation of $\fg$, that is $\K^2$ equipped with the left action of $\fg = \sl_2(\K)$.
\end{notation*}

\begin{notation*}
For $k \geq 1$ an integer, let $\Sym^k \Nat \fg$\inmargin{$\Sym^k \Nat \fg$} denote the $k\th$ symmetric power of $\Nat \fg$.
\end{notation*}

We do not wish to go into tensor algebra, and will more conveniently handle $\Sym^k \Nat \fg$ as follows.

\begin{fact*}[{\cite[\S II.7]{Humphreys}}]
Let $\K$ be a field of characteristic $0$ or $\geq k+1$. Then:
\begin{itemize}
\item
$S_k = \Sym^k \Nat \fg$ is isomorphic to $\K[Y^k, XY^{k-1}, \dots, X^k]$ as a $\K\fg$-module, where $x$ acts as $X \frac{\partial}{\partial Y}$ and $y$ as $Y \frac{\partial}{\partial X}$;
\item
$S_k$ is an irreducible $\K\fg$-module; it remains irreducible as a $\fg$-module;
\item
$h = [x, y]$ acts on $S_k$ as $X\frac{\partial}{\partial X} - Y\frac{\partial}{\partial Y}$;
\item
$\K\cdot X^{k-i} Y^i = E_{k-2i}(S_k)$;
\item
the length of $S_k$ is $k+1$, meaning that $\fu^{k+1}\cdot S_k = 0$ and $\fu^k \cdot S_k \neq 0$;
\item
The Casimir operator $c_1$ acts on $S_k$ as multiplication by $k(k+2)$. In particular in characteristic $0$ or $\geq k+3$, $c_1$ induces a bijection of $S_k$.
\end{itemize}
\end{fact*}

\section{Length}\label{S:length}

This section contains two minor results on the notion of length as defined in \S\ref{S:cadre}. They are fairly straightforward and so are their proofs.

\begin{variation}\label{v:xlambdan=0lengthn}
Let $n \geq 2$ be an integer and $\K$ be a field of characteristic $0$ or $\geq n+1$.
Let $\fg = \sl_2(\K)$, $\fb \leq \fg$ be a Borel subring and $\fu = \fb'$ be its radical. Let $V$ be a $\fu$-module.
Suppose that for all $\lambda \in \K$, $x_\lambda^n = 0$ in $\End V$. Then $V$ has $\fu$-length at most $n$, meaning that $\fu^n \cdot V = 0$.
\end{variation}
\begin{proof}
This is a simpler analog of Variation n$^\circ$6 \cite{TV-I}: only the end of the argument need be reproduced as the induction on the ``weights'' of monomials is not necessary. Indeed $x_{\lambda + \mu} =  x_\lambda + x_\mu$ whence immediately:
\[0 = \sum_{j = 1}^{n-1} \binom{n}{j} x_\lambda^j x_\mu^{n-j}\]

One then replaces $\mu$ by $i \mu$ for $i = 1\dots n-1$; this yields the same linear $(n-1)$ by $(n-1)$ system as in Variation n$^\circ$6 \cite{TV-I}. Hence $d x_\lambda^{n-1} x_\mu = 0$ where $d$ is the determinant of the system; all prime divisors of $d$ divide $n!$. In particular replacing $\mu$ by $\frac{\mu}{d}$ in $\K$, one finds $x_\lambda^{n-1} x_\mu = 0$ in $\End V$. Since $\fu$ acts on $\im x_\mu$ one may use induction on $n$.
\end{proof}

%
%

Next comes an easy generalization of Variation n$^\circ$9 \cite{TV-I}. (There seems to be no parallel argument in the case of $\SL_2(\K)$; the quadratic setting painfully dealt with in Variation n$^\circ$7 \cite{TV-I} actually required a full $\SL_2(\K)$-module.)

\begin{variation}\label{v:xn=0l=0}
Let $n \geq 2$ be an integer and $\K$ be a field of characteristic $0$ or $\geq n+1$.
Let $\fg = \sl_2(\K)$ and $\fb \leq \fg$ be a Borel subring. Let $V$ be a $\fb$-module. Suppose that $x^n \cdot V = 0$. Then $V$ has $\fu$-length at most $n$, meaning that $\fu^n \cdot V = 0$.
\end{variation}
\begin{proof}
We go to $\End V$. Let us prove by induction on $i = 0 \dots n$:
\[\forall (\lambda_1, \dots, \lambda_i) \in \K^i, \ x^{n-i} x_{\lambda_i} \dots x_{\lambda_1} = 0\]
\begin{itemize}
\item
This holds of $i = 0$.
\item
Suppose that the result holds of fixed $i < n$; let $(\lambda_1, \dots, \lambda_i, \lambda_{i+1}) \in \K^{i+1}$.

We show by induction on $j = 0 \dots i$:
\[x^{n-i} x_{\lambda_i} \dots x_{\lambda_{j+1}} h_{\lambda_{i+1}} x_{\lambda_j} \dots x_{\lambda_1} = 0\]
\begin{itemize}
\item
This holds of $j = 0$ by assumption on $i$.
\item
Suppose that the result holds of fixed $j$. Then:
\begin{align*}
& \quad x^{n-i} x_{\lambda_i} \dots x_{\lambda_{j+2}} h_{\lambda_{i+1}} x_{\lambda_{j+1}} \dots x_{\lambda_1}\\
& = x^{n-i} x_{\lambda_i} \dots x_{\lambda_{j+2}} ([h_{\lambda_{i+1}} , x_{\lambda_{j+1}} ] + x_{\lambda_{j+1}}h_{\lambda_{i+1}}) x_{\lambda_j} \dots x_{\lambda_1}\\
& = \quad 2 x^{n-i} x_{\lambda_i} \dots x_{\lambda_{j+2}} x_{\lambda_{i+1}\cdot \lambda_{j+1}} x_{\lambda_j} \dots x_{\lambda_1}\\& \quad + x^{n-i} x_{\lambda_i} \dots x_{\lambda_{j+1}} h_{\lambda_{i+1}} x_{\lambda_j} \dots x_{\lambda_1}\\
& = 0
\end{align*}
by assumption on $j$ and $i$ (the latter applied with $\lambda'_{j+1} = \lambda_{i+1} \cdot \lambda_{j+1}$). This concludes the induction on $j$.
\end{itemize}
With $j = i$, one gets:
\[x^{n-i} h_{\lambda_{i+1}} x_{\lambda_i} \dots x_{\lambda_1} = 0\]
Let us now prove by induction on $k = 0 \dots n-i$:
\[x^{n-(i+k)} h_{\lambda_{i+1}} x^k x_{\lambda_i} \dots x_{\lambda_1} = 2k x^{n-(i+1)} x_{\lambda_{i+1}} \dots x_{\lambda_1}\]
\begin{itemize}
\item
This holds of $k = 0$ by what we have just shown.
\item
Suppose that the result holds of fixed $k$. Then:
\begin{align*}
& \quad x^{n-(i + k+ 1)} h_{\lambda_{i+1}} x^{k+1} x_{\lambda_i} \dots x_{\lambda_1} \\
& = x^{n-(i + k+ 1)} ([h_{\lambda_{i+1}}, x] + x h_{\lambda_{i+1}}) x^{k} x_{\lambda_i} \dots x_{\lambda_1}\\
& = 2 x^{n-(i+k+1)} x_{\lambda_{i+1}} x^k x_{\lambda_i} \dots x_{\lambda_1} + x^{n-(i+k)} h_{\lambda_{i+1}} x^k x_{\lambda_i} \dots x_{\lambda_1}\\
& = 2 x^{n-(i+1)}x_{\lambda_{i+1}} \dots x_{\lambda_1} + 2k x^{n-(i+1)}x_{\lambda_{i+1}} \dots x_{\lambda_1}
\end{align*}
This concludes the induction on $k$.
\end{itemize}
With $k = n-i$ one gets:
\[h_{\lambda_{i+1}} x^{n-i} x_{\lambda_i} \dots x_{\lambda_1} = 2(n-i) x^{n-(i+1)} x_{\lambda_{i+1}}\dots x_{\lambda_i}\]
but the left-hand side is zero by assumption. If instead of $\lambda_{i+1}$ we had started with $\frac{\lambda_{i+1}}{2(n-i)}$, which is legitimate by assumption on the characteristic of $\K$, we would have obtained:
\[x^{n-(i+1)} x_{\lambda_{i+1}} \dots x_{\lambda_1} = 0\]
This concludes the induction on $i$.
\end{itemize}
With $i = n$, one has the desired statement.
\end{proof}

\begin{remark*}
One cannot use induction on $n$ via $\im x$ since $\im x$ may fail to be $\ft = \{h_\lambda: \lambda \in \K\}$-invariant; such a configuration will be met in the example illustrating the following.
\end{remark*}

\begin{remark*}
The mere existence of a product $x_{\lambda_1}\dots x_{\lambda_n}$ which is zero in $\End V$ does not suffice to force the length to be at most $n$.
Take indeed $\K = \C$, $\fg = \sl_2(\C)$, and let $\varphi$ stand for complex conjugation. Also let $V = \Nat \fg \simeq \C^2$, $V' = {}^\varphi V$ (a copy ``twisted'' by the field automorphism), and $W = V \otimes V'$. One sees that $W$ has no $\C\fg$-submodule other than $\{0\}$ and $W$.

Let $(e_1, e_2)$ be the standard basis of $\C^2$; one has $x \cdot e_1 = 0$ and $x\cdot e_2 = e_1$. Write for simplicity $e_{i, j} = e_i \otimes e_j$. One finds:
\[\left\{\begin{array}{rcl}
x_\lambda \cdot e_{1,1} & = & 0\\
x_\lambda \cdot e_{2,1} & = & \lambda e_{1,1}\\
x_\lambda \cdot e_{1,2} & = & \varphi( \lambda) e_{1,1}\\
x_\lambda \cdot e_{2,2} & = & \lambda e_{1,2} + \varphi (\lambda) e_{2,1}
\end{array}\right.\]
so that $x_\lambda x_\mu \cdot e_{2,2} = (\lambda \varphi (\mu) + \mu \varphi (\lambda)) e_{1,1}$.
Clearly $x_1^2 \neq 0$ and yet $x_1 x_i = 0$ (where $i$ stands for a root of $-1$).

One may object that $W$ though simple as a $\C\fg$-module, is not as a $\fg$-module; we then go down to $W_0 = \R e_{1, 1} \oplus \{\lambda e_{1,2} + \varphi(\lambda) e_{2,1}: \lambda \in \C\} \oplus \R e_{2,2}$, which as a $\fg$-module is simple; one has $x_1^2 \neq 0$ in $\End W_0$.
\end{remark*}

\begin{remark*}
By Variation \ref{v:xn=0l=0} the $\fu$-length is therefore the nilpotence order of $x$ in $\End V$; one may wonder whether it is the nilpotence order of $y$ as well.
(One should not expect this in full generality: after Variation n$^\circ$12 \cite{TV-I} we saw that it can be achieved in characteristic $3$ that $x^2 = 0 \neq y^2$.)

Here is an unsatisfactory argument in characteristic zero.
\begin{quote}
We go to the enveloping algebra $\fA$. Then $\Ann_{\fA}(V)$ is a two-sided ideal containing $x^n$. But $\SL_2(\K)$ acts on $\fA$ and normalizes every two-sided ideal by \cite[Proposition 2.4.17]{Dixmier}; since the Weyl group exchanges $x$ and $y$ one has $y^n \in \Ann_\fA(V)$ as well, whence $y^n = 0$ in $\End V$.
\end{quote}

The argument is not quite satisfactory: we have been using the $\K$-algebra $\fA$. It is a fact that every $\K\fg$-module is an $\fA$-module but all we had was a mere $\fg$-module; turning it into a $\K\fg$-module is precisely the core of the matter.

More prosaically, a crude computation will show that in characteristic $\geq 2n+1$ one does have $x^n = 0 \Rightarrow y^n = 0$: we shall see this while proving Variation \ref{v:oplusSymNatsl2K1}.
\end{remark*}

\section{Combinatorial skeleton}\label{S:squelette}

In this section we focus on $\sl_2(\K_1)$-modules of finite length, with $\K_1$ a prime field. If the characteristic is $0$ or large enough, Variation \ref{v:oplusSymNatsl2K1} of \S\ref{s:pgeq2n+1} gives a complete description. But some other objects appear if one tries to lower the characteristic too much (\S\ref{s:pathologies}). Provided one assumes that $y$ has the same order as $x$, the monsters vanish (Variation \ref{v:oplusSymNatsl2K1:lowchar}, \S\ref{s:xn=yn=0}).

The author cannot believe that the results of this section are new, but found no evidence. We shall give purely computational arguments without going to the algebraic closure, which was another possible direction.

\subsection{Large Enough Characteristic}\label{s:pgeq2n+1}

\begin{variation}\label{v:oplusSymNatsl2K1}
Let $n \geq 2$ be an integer and $\K_1$ be a \emph{prime} field of characteristic $0$ or $\geq 2n+1$.
Let $\fg_1 = \sl_2(\K_1)$ and $V$ be a $\fg_1$-module. If the characteristic of $\K$ is $0$ one requires $V$ to be torsion-free. Suppose that $x^n = 0$ in $\End V$. 

Then $V = \Ann_V(\fg_1) \oplus \fg_1 \cdot V$, and $\fg_1\cdot V$ is a $\K_1$-vector space with $\fg_1 \cdot V \simeq \oplus_{k = 1}^{n-1} \oplus_{I_k} \Sym^k \Nat \fg_1$ as $\K_1\fg_1$-modules.
\end{variation}
\begin{proof}
Induction on $n$. When $n = 2$ this is Variation n$^\circ$12 \cite{TV-I}.

All along $c_1 = 2xy+2yx+h^2$ will be the Casimir operator; the action of $c_1$ commutes with the action of $\fg_1$ on $V$. In characteristic $0$ this holds only since we assumed $V$ to be torsion-free.

\begin{step}[see Variation n$^\circ$3 \cite{TV-I}]\label{v:oplusSymNatsl2K1:st:W}
We may assume $V = \fg_1 \cdot V$ and $\Ann_V(\fg_1) = 0$.
\end{step}
\begin{proofclaim}
Let $W = \fg_1\cdot V$ and $\overline{W} = W/\Ann_W (\fg_1)$. Let $\overline{\cdot}$ stand for projection modulo $\Ann_W \fg_1$. By perfectness of $\fg_1$ one has $\overline{W} = \fg_1\cdot \overline{W}$ and $\Ann_{\overline{W}}\fg_1 = 0$. In particular if $\K_1 = \Q$ then $\overline{W}$ is torsion-free. Suppose that the result holds of $\overline{W}$; let us prove it for $V$: suppose that $\overline{W}$ is a $\K_1$-vector space satisfying $\overline{W} \simeq \oplus_{k = 1}^{n-1} \oplus_{I_k} \Sym^{n-1} \Nat \fg_1$ as $\K_1\fg_1$-modules.

One then sees that $c_1$ is a bijection of $\overline{W}$. We claim the following:
\begin{itemize}
\item
$W = c_1\cdot W + \Ann_W \fg_1$. For take $w \in W$. Since $c_1$ is surjective onto $\overline{W}$ there exists $w' \in W$ with $\overline{w} = c_1 \cdot \overline{w'}$.
\item
$c_1 \cdot W = W$. Let us apply $\fg_1$ to the previous equality, bearing in mind perfectness of $\fg_1$ and centralness of $c_1$. One finds $W = \fg_1 \cdot W = \fg_1 c_1 \cdot W + \fg_1 \cdot \Ann_W \fg_1 = c_1 \fg_1 \cdot W = c_1\cdot W$.
\item
$W \cap \ker c_1 = 0$. For take $w \in W \cap \ker c_1$. Then by the previous claim there exists $w' \in W$ with $w = c_1 \cdot w'$. Modulo $\Ann_W \fg_1$ one has $0 = c_1 \cdot \overline{w} = c_1^2 \cdot \overline{w'}$. By injectivity of the Casimir operator on $\overline{W}$ it follows $\overline{w'} = 0$, whence $w' \in \Ann_V \fg_1 \leq \ker c_1$ and $w = c_1 \cdot w' = 0$.
\item
$\Ann_V \fg_1 = \ker c_1$. One inclusion is obvious and was just used; if conversely $k \in \ker c_1$ then $\fg_1\cdot \<c_1\cdot k\> = 0 = c_1\cdot \<\fg_1 \cdot k\>$ so $\<\fg_1\cdot k\> \leq W \cap \ker c_1 = 0$.
\item
$V = \Ann_V \fg_1 \oplus W$. The sum is direct indeed as we just saw. Moreover if $v \in V$ then there exists $w \in W$ with $c_1 \cdot v = c_1 \cdot w$; in particular $V \leq W + \ker c_1 = W \oplus \ker c_1 = W \oplus \Ann_V \fg_1$.
\end{itemize}
$V$ therefore has the desired structure.
\end{proofclaim}

We now suppose $V = \fg_1 \cdot V$ and $\Ann_V \fg_1 = 0$.
It follows that $V \simeq (\fg_1\cdot V)/(\Ann_{\fg_1\cdot V}\fg_1)$ is a $\K_1$-vector space.

\begin{step}\label{v:oplusSymNatsl2K1:st:h}
In $\End V$, $(h-n+1) (h-n+2) \dots (h+n-1) = 0$.
\end{step}
\begin{proofclaim}
Remember that in the enveloping ring, for $i, j \geq 1$:
\[
y^jx^i = \sum_{k = 0}^{\min(i, j)} \left[(-1)^k k!\ \binom{i}{k} \binom{j}{k} \left(\prod_{\ell = 0}^{k-1}(h-i+j+\ell)\right) x^{i-k}y^{j-k}\right]
\]

In the subring of $\End V$ generated by the image of $\fg_1$ one has $x^n = 0$; the formula becomes with $i=n$ and $j \leq n$:
\[\sum_{k=1}^j (-1)^k k!\ \binom{n}{k} \binom{j}{k} \left(\prod_{\ell = 0}^{k-1} (h-n+j+\ell)\right) x^{n-k} y^{j-k} = 0\qquad (F_j)\]
Let us prove by induction on $j = 0 \dots n$:
\[(h-n+1) (h-n+2) \dots (h-n+2j-1) x^{n-j} = 0\]
When $j = 0$ the (ascending) product is empty: our claim holds by assumption on $x$. Suppose that the result holds of fixed $j$ and let us prove it for $j+1 \leq n$.
Consider formula $(F_{j+1})$ multiplied on the left by $(h-n+1)\dots (h-n+j)$. One gets:
\[\sum_{k = 1}^{j+1} (-1)^k k!\ \binom{n}{k} \binom{j+1}{k} \pi_k x^{n-k} y^{j+1-k} = 0\]
where:
\begin{align*}
\pi_k & = (h-n+1)\dots (h-n+j)\cdot (h-n+j+1) \dots (h-n+j+k)\\
& = (h-n+1) \dots (h-n + j + k)
\end{align*}
Since $j + k \geq 2 k -1$ the term with index $k$ contains $(h-n+1) \dots (h-n+2k-1) x^{n-k}$, which by induction is zero while $k \leq j$. So only remains the term with index $k = j+1$ namely:
\[(-1)^{j+1} (j+1)!\ \binom{n}{j+1} (h-n+1) \dots (h-n+2j+1) x^{n-(j+1)} = 0\]
By $n!$-torsion-freeness of $V$ we may remove the coefficients and complete the induction.
When $j = n$ one finds $(h-n+1) (h-n+2) \dots (h+n-1) = 0$.
\end{proofclaim}

\begin{step}\label{v:oplusSymNatsl2K1:st:Ei}
$V = \oplus_{j = 1-n}^{n-1} E_j$.
\end{step}
\begin{proofclaim}
Let us first observe that the sum $\oplus_{j = 1-n}^{n-1} E_j$ is direct indeed by $(2n-2)!$-torsion-freeness of $V$. For the same reason the monomials $X-j$ are pairwise coprime in $\K_1[X]$ for $j = 1 - n, \dots, n-1$. Since their product annihilates $h$ in $\End V$ one has $V = \oplus_{j = 1-n}^{n-1} \ker (h-j) = \oplus_{j = 1-n}^{n-1} E_j$.
\end{proofclaim}

Since $n-1+2 = n+1$ and $n-2 + 2 = n$ are not congruent to any $j \in \{1-n, \dots, n-1\}$ the operator $x$ annihilates $E_{n-1}$ and $E_{n-2}$. Similarly $y$ annihilates $E_{1-n}$ and $E_{2-n}$.

\begin{localremark*}
It is now clear that $y^n \cdot V = 0$.
\end{localremark*}

\begin{localnotation}\label{v:oplusSymNatsl2K1:n:VbotVtop}
Let $V_\bot = \im(c_1 - n^2 +1)$ and $V_\top = \ker (c_1 - n^2 + 1)$.
\end{localnotation}

\begin{step}\label{v:oplusSymNatsl2K1:st:Vbot}
$V_\bot$ is a $\fg_1$-submodule isomorphic to $\oplus_{k = 1}^{n-2}\oplus_{I_k}\Sym^k \Nat \fg_1$.
\end{step}
\begin{proofclaim}
$V_\bot$ is clearly $\fg_1$-invariant.
But by Step \ref{v:oplusSymNatsl2K1:st:Ei} or the proof of Step \ref{v:oplusSymNatsl2K1:st:h} one has in $\End V$ the identity $hx^{n-1} = (n-1) x^{n-1}$. Hence always in $\End V$:
\begin{align*}
x^{n-1} (c_1-n^2 + 1) & = x^{n-1} (2xy+2yx+h^2 - n^2 + 1)\\
& = 2(x^{n-1}y)x +(h+2-2n)^2 x^{n-1} - (n^2-1)x^{n-1}\\
& = 2 (yx^{n-1} +(n-1)(h+2-n)x^{n-2})x\\
& \quad + (1-n)^2 x^{n-1} + (1-n^2)x^{n-1}\\
& = 2(n-1) (h+2-n) x^{n-1} + 2 (1-n) x^{n-1}\\
& = 0
\end{align*}
It follows that $x^{n-1}$ annihilates $\im (c_1-n^2 + 1) = V_\bot$ and one may apply induction. Since $\Ann_{V_\bot}\fg_1 \leq \Ann_V \fg_1 = 0$ there remains only $V_\bot = \fg_1\cdot V_\bot \simeq \oplus_{k = 1}^{n-2}\oplus_{I_k}\Sym^k \Nat \fg_1$.
\end{proofclaim}

\begin{step}\label{v:oplusSymNatsl2K1:st:V=Vbot+Vtop}
We may assume $V = V_\top$.
\end{step}
\begin{proofclaim}
We claim that $V = V_\bot \oplus V_\top$. The way the Casimir operator acts on each $\Sym^k \Nat \fg_1$ is known: like multiplication by $k(k+2)$. But for $k = 1, \dots, n-2$, $k(k+2) \neq n^2 - 1$ in $\K_1$ by assumption on the characteristic. It follows that $(c_1 - n^2 + 1)$ induces a bijection of $V_\bot$. As a consequence $V_\bot \cap V_\top = V_\bot \cap \ker (c_1 - n^2 + 1) = 0$. Moreover for all $v \in V$ there exists $v' \in V_\bot$ with $(c_1 - n^2 + 1)\cdot v = (c_1 - n^2 + 1)\cdot v'$, whence $V = V_\bot + \ker (c_1 - n^2 + 1) = V_\bot + V_\top = V_\bot \oplus V_\top$.
If the result were proved for $V_\top$ it would therefore follow for $V$.
\end{proofclaim}

From now on we suppose $V = V_\top$; in particular $c_1 - n^2 + 1$ annihilates $V$.

\begin{step}\label{v:oplusSymNatsl2K1:st:Vtopsepare}
$\ker x = E_{n-1}$.
\end{step}
\begin{proofclaim}
We claim that $x$ is injective on $\oplus_{j = 1-n}^{n-2} E_j$. For if $a \in E_j$ with $j \in \{1-n, \dots, n-2\}$ satisfies $x \cdot a = 0$, then
\[
(n^2 - 1) a = c_1 \cdot a = (2xy+2yx+h^2) \cdot a = (2h + h^2) \cdot a = j(j+2) a
\]
so either $a = 0$ or $n^2 = j(j+2)+1 = (j+1)^2$. But the latter equation solves into $j = \pm n -1$ which is not the case (even in characteristic $p \geq 2n+1$).
\end{proofclaim}

\begin{step}\label{v:oplusSymNatsl2K1:st:Vtop}
$V$ is isomorphic to $\oplus_{I_{n-1}} \Sym^{n-1} \Nat\fg_1$.
\end{step}
\begin{proofclaim}
%
We claim that for all $i = 1 \dots n$, $E_{n-2i} = 0$. At $i = 1$ this is because $E_{n-2} \leq \ker x = E_{n-1}$ by Step \ref{v:oplusSymNatsl2K1:st:Vtopsepare}. If this is known at $i$, then $x \cdot E_{n-2(i+1)} \leq E_{n-2i}$ whence $E_{n-2(i+1)} \leq \ker x = E_{n-1}$.

On the other hand observe that for all $i = 1 \dots n$:
$yx_{|E_{n + 1 - 2 i}} = (i-1)(n + 1 - i)$ and $xy_{|E_{n + 1 - 2 i}} = i(n - i)$.
This is actually obvious since $c_1 = 4 yx + h^2 + 2 h = 4 xy + h^2 - 2h$ is constant and equals multiplication by $n^2 - 1$.

It is therefore now clear that for all $a_{n-1} \in E_{n-1}\setminus\{0\}$, $\<\fg_1\cdot a_{n-1}\>$ is a $\K_1$-vector space isomorphic to $\Sym^{n-1}\Nat \fg_1$ as a $\K_1\fg_1$-module; if in particular $b \in \<\fg_1\cdot a_{n-1}\>\setminus \{0\}$ then $\<\fg_1\cdot a_{n-1}\> = \<\fg_1\cdot b\>$. Let $M \leq V$ be a maximal direct sum of such spaces. Then $M$ has the desired structure, and our computations show $V = \oplus_{j = 1-n}^{n-1} E_j \leq \fg_1 \cdot E_{n-1} \leq M$.
\end{proofclaim}
This finishes the proof.
\end{proof}

\subsection{A Digression: Pathologies in Low Characteristic}\label{s:pathologies}

As in Variation n$^\circ$12 \cite{TV-I} the characteristic must be $0$ or $\geq 2n+1$ in order to prove Variation \ref{v:oplusSymNatsl2K1}. In this section we suppose the characteristic to lie between $n$ and $2n$. We shall construct counterexamples to Variation \ref{v:oplusSymNatsl2K1} and remove them later in \S\ref{s:xn=yn=0} under the extra assumption that $y$ has the same order as $x$ in $\End V$.

The construction generalizes the one given in characteristic $3$ at the end of \cite{TV-I}.
Let $n \geq 2$ be an integer and $p$ be a prime number with $n < p < 2n$; let $m$ be such that $n + m = p$. Observe that if $i \in \{1, \dots, n \}$ and $j \in \{1, \dots, m\}$, then $n+1-2i$ and $m+1-2j$ are never congruent modulo $p$. Hence modulo $p$, the $n+1-2i$'s and $m+1-2j$'s are all distinct, and their global number is $p$; there are $n$ of the former kind and $m$ of the latter.

\begin{construction*}
Let $V_1$ and $V_2$ be two vector spaces over $\F_p$. Let $\alpha: V_1 \to V_2$ and $\beta: V_2 \to V_1$ be two linear maps. Define a $\fg_1$-module $\bS_{\alpha, \beta}$ as follows.

For each $j \in \{1, \dots, m\}$ let $E_{m+1-2j}$ be a copy of $V_1$ whose elements we shall denote $e_{m+1-2j, v_1}$ for $v_1 \in V_1$. For each $i \in \{1, \dots, n\}$ let $E_{n+1-2i} = \{e_{n+1-2i, v_2} : v_2 \in V_2\}$ be a copy of $V_2$.

The underlying vector space of $\bS_{\alpha, \beta}$ is
\[\left(\mathop\oplus_{j = 1}^m E_{m+1-2j} \right)\bigoplus \left(\mathop\oplus_{i = 1}^n E_{n+1-2i}\right)\] Now define an action of $\fg_1 = \sl_2(\F_p)$ on $\bS_{\alpha, \beta}$ by:
\begin{align*}
h \cdot e_{m+1-2j, v_1} & = (m+1-2j) \cdot e_{m+1-2j, v_1} & \mbox{ if } 1 \leq j \leq m\\
x \cdot e_{m+1-2j, v_1} & = (j-1) e_{m+1-2(j-1), v_1} & \mbox{ if } 1 \leq j \leq m\\
y \cdot e_{m+1-2j, v_1} & = (m-j) e_{m+1-2(j+1), v_1} & \mbox{ if } 1 \leq j < m\\
y \cdot e_{1-m, v_1} & = e_{n-1, \alpha(v_1)}\\
h \cdot e_{n+1-2i, v_2} & = (n+1-2i) \cdot e_{n+1-2i, v_2} & \mbox{ if } 1 \leq i \leq n\\
x \cdot e_{n+1-2i, v_2} & = (i-1) e_{n+1-2(i-1), v_2} & \mbox{ if } 1 \leq i \leq n\\
y \cdot e_{n+1-2i, v_2} & = (n-i) e_{n+1-2(i+1), v_2} & \mbox{ if } 1 \leq i < n\\
y \cdot e_{1-n, v_2} & = e_{m-1, \beta(v_2)}
\end{align*}
\end{construction*}

Note that by construction, $x$ annihilates $E_{n-1}$ and $E_{m-1}$.

\[\begin{tikzpicture}[scale=.3]
  \node (1-n) at (0,0) {$E_{1-n}$};
  \node (n-1) at (16,0) {$E_{n-1}$};
  \node (m-1) at (2,4) {$E_{m-1}$};
  \node (1-m) at (14,4) {$E_{1-m}$};
  \node (bdn) at (4,0) {};
  \node (edn) at (12,0) {};
  \node (bdm) at (10,4) {};
  \node (edm) at (6,4) {};
  \draw[->] (1-n) to node {$x$} (bdn);
  \draw[->] (edn) to node {$x$} (n-1);
  \draw[loosely dotted] (bdn) to (edn);
\draw[->] (1-m) to node {$x$} (bdm);
  \draw[->] (edm) to node {$x$} (m-1);
  \draw[loosely dotted] (bdm) to (edm);
  \draw[->, bend left] (bdn) to node {$y$} (1-n.352);
  \draw[->, bend left] (n-1.188) to node {$y$} (edn);
  \draw[->, bend left] (bdm) to node {$y$} (1-m.172);
  \draw[->, bend left] (m-1.8) to node {$y$} (edm);
\draw[->, bend left] (1-m) to node {$\alpha$} (n-1);
\draw[->, bend left] (1-n) to node {$\beta$} (m-1);
\draw (22,-0.5) node {$\left\}\begin{array}{c}\bigoplus\Sym^{n-1}\\\Nat\fg_1\end{array}\right.$};
\draw (-4,4.5) node {$\left.\begin{array}{c}\bigoplus\Sym^{m-1}\\\Nat\fg_1\end{array}\right\{$};
\end{tikzpicture}
\]

\begin{observation*}
$\bS_{\alpha, \beta}$ is a $\fg_1$-module annihilated by $x^n$.
\end{observation*}
\begin{proofclaim}
It suffices to prove that the defining relations of $\fg_1$ are satisfied at every vector $e_{m+1-2j, v_1}$; the $e_{n+1-2i, v_2}$'s are treated similarly.

At $e_{m+1-2j, v_1}$ with $1 < j < m$ there is nothing to prove since everything is locally as in $\Sym^{n-1}\Nat\fg_1$; since $x$ annihilates $e_{m-1, v_1}$, this also holds at $e_{m-1, v_1}$.
Let us now consider a vector $e_{1-m, v_1}$.
One checks:
\begin{align*}
x \cdot (y \cdot e_{1-m, v_1}) - y \cdot (x \cdot e_{1-m, v_1}) & = x \cdot e_{n-1, \alpha(v_1)} - y \cdot (m-1) e_{3-m, v_1}\\
& = (1-m) e_{1-m, v_1}\\
& = h \cdot e_{1-m, v_1}
\end{align*}
then $h\cdot (x\cdot e_{1-m, v_1}) - x\cdot (h\cdot e_{1-m, v_1}) = 0 = 2x \cdot e_{1-m, v_1}$, and finally:
\begin{align*}
h\cdot (y\cdot e_{1-m, v_1}) - y\cdot (h\cdot e_{1-m, v_1}) & = h \cdot e_{n-1, \alpha(v_1)} + (m-1) y \cdot e_{1-m, v_1}\\
& = (n-1) e_{n-1, \alpha(v_1)} + (m-1) e_{n-1, \alpha(v_1)}\\
& = (p - 2) e_{n-1, \alpha(v_1)}\\
& = - 2 y \cdot e_{1-m, v_1}
\end{align*}
By construction, $x^n \cdot \bS_{\alpha, \beta} = 0$.
\end{proofclaim}


Our construction could a priori depend on bases we chose for $V_1$ and $V_2$; it is up to isomorphism not the case.

\begin{observation*}
$\bS_{\alpha, \beta}$ and $\bS_{\alpha', \beta'}$ are isomorphic iff the pairs $(\alpha, \beta)$ and $(\alpha', \beta')$ are equivalent, that is iff there exist linear isomorphisms $u_1: V_1 \simeq V'_1$ and $u_2: V_2 \simeq V'_2$ with $\alpha' = u_2 \alpha u_1^{-1}$ and $\beta' = u_1 \beta u_2^{-1}$.
\end{observation*}
\begin{proofclaim}
If the pairs $(\alpha, \beta)$ and $(\alpha', \beta')$ are equivalent, an isomorphism of $\fg_1$-modules $\bS_{\alpha, \beta} \simeq \bS_{\alpha', \beta'}$ is easily constructed by setting $f(e_{m+1-2j, v_1}) = e'_{m+1-2j, u_1(v_1)}$ with obvious notations and similarly on the other row.

For the converse suppose that there is such an isomorphism $f : \bS_{\alpha, \beta} \simeq \bS_{\alpha', \beta'}$.
Let $V_1 = E_{1-m}(\bS_{\alpha, \beta})$ and $V_2 = E_{1-n}(\bS_{\alpha, \beta})$; these are $\F_p$-vector spaces. One then retrieves $\alpha(v_1) = \frac{1}{(n-1)!} y^n \cdot v_1$ and $\beta(v_2) = \frac{1}{(m-1)!}y^m \cdot v_2$ (which do induce $\bS_{\alpha, \beta}$). Proceed similarly on $\bS_{\alpha', \beta'}$.

Let $u_1(v_1) = f(v_1)$ and $u_2(v_2) = f(v_2)$. Since $f$ is an isomorphism they are linear isomorphisms between $V_1$ and $V'_1$, resp. $V_2$ and $V'_2$. Now for all $v_1 \in V_1$, since $f$ is an isomorphism of $\fg_1$-modules,
\begin{align*}
u_2 \circ \alpha (v_1) & = u_2 \left(\frac{1}{(n-1)!} y^n \cdot v_1\right)\\
& = \frac{1}{(n-1)!} f(y^n \cdot v_1)\\
& = \frac{1}{(n-1)!} y^n \cdot f(v_1)\\
& = \alpha'(f(v_1))\\
& = \alpha'\circ u_1(v_1)
\end{align*}
A similar verification can be carried on $V_2$, proving that $u_1$ and $u_2$ define a equivalence of $(\alpha, \beta)$ and $(\alpha', \beta')$.
\end{proofclaim}

\begin{observation*}
$\bS_{\alpha, \beta}$ is \emph{non-}simple iff there are subspaces $W_1 \leq V_1$ and $W_2 \leq V_2$ \emph{not both zero} such that $\alpha$ maps $W_1$ to $W_2$ and $\beta$ maps $W_2$ to $W_1$.
\end{observation*}
\begin{proofclaim}
We give a correspondence between submodules of $\bS_{\alpha, \beta}$ and pairs $(W_1, W_2)$ as in the statement. One direction is clear: if such a pair $(W_1, W_2)$ is given, a $\fg_1$-submodule is readily defined.

So let $W \leq V$ be a $\fg_1$-submodule. Set $W_1 = \{v_1 \in V_1 : e_{1-m, v_1} \in W\}$ and $W_2 = \{v_2 \in V_2 : e_{1-n, v_2} \in W\}$.
We claim that $\alpha$ maps $W_1$ to $W_2$, and that $\beta$ maps $W_2$ to $W_1$. It suffices to prove it for $\alpha$.
Take indeed $w_1 \in W_1$. Then by construction $e_{1-m, w_1} \in W$ whence $y \cdot e_{1-m, w_1} = e_{n-1, \alpha(w_1)} \in W$.
Applying $y^{n-1}$, one finds up to multiplication by $(n-1)!$ which is coprime to $p$: $e_{1-n, \alpha(w_1)} \in W$, so by definition, $\alpha(w_1) \in W_2$.
\end{proofclaim}

\begin{observation*}
$\bS_{\alpha, \beta}$ is simple iff $\alpha$ and $\beta$ are isomorphisms and $\beta\circ\alpha$ is irreducible (as an automorphism of $V_1$).
\end{observation*}
\begin{proofclaim}
Suppose that $\bS_{\alpha, \beta}$ is simple. Take $f\in V_1 \setminus\{0\}$. Consider the sequences $f_r = (\beta\alpha)^r(f)$ and $g_r = \alpha(f_r)$. These sequences span subspaces $W_1 \leq V_1$ and $W_2 \leq V_2$ mapped one to another by $\alpha$ and $\beta$. By simplicity $W_1 = V_1$ and $W_2 = V_2$. Hence $\alpha$ is a bijection. A similar argument holds for $\beta$.
Now let $W_1$ be invariant under $\beta\circ\alpha$ and set $W_2 = \alpha(W_1)$. Then $\alpha$ maps $W_1$ to $W_2$ and $\beta$ maps $W_2$ to $W_1$ so by simplicity $W_1 = V_1$ or $0$.

Suppose that $\alpha$ and $\beta$ are isomorphisms such that $\beta\circ\alpha$ is irreducible. If $W_1 \leq V_1$ and $W_2\leq V_2$ are mapped one to another by $\alpha$ and $\beta$ then $\beta\circ \alpha (W_1) \leq W_1$ so $W_1 = 0$ or $V_1$. In the former case $W_2 = 0$ by injectivity of $\beta$, in the latter $W_2 = V_2$ by surjectivity of $\alpha$.
\end{proofclaim}

Suppose in particular that $\bS_{\alpha, \beta}$ is simple. Then for any $f \in V_1\setminus\{0\}$, the sequence $(f_r)_{r\geq 1}$ as above spans $V_1$ whence a linear relation $f_0 = \sum_{r = 1}^d k_r f_r$. It follows that $f_d$ lies in the span $W_1$ of $(f_r)_{0 \leq r \leq d-1}$, so $V = W_1$ is finite-dimensional. Moreover the characteristic polynomial of $\beta\circ\alpha$ is irreducible over $\F_p[X]$.

\begin{observation*}
Let $n \geq 2$ be an integer and $\K_1$ be the field $\F_p$ with $n < p < 2n$.
Let $\fg_1 = \sl_2(\K_1)$ and $V$ be a \emph{simple} $\fg_1$-module. Suppose that $x^n = 0$ in $\End V$.

Then $V$ is some $\bS_{\alpha, \beta}$.
\end{observation*}
\begin{proofclaim}
We may suppose $n$ minimal such that $x^n = 0$. (The reader will observe that had we wished to be fully rigorous we should have written $\bS_{\alpha, \beta}^n$ throughout.)

By simplicity $\Ann_V \fg_1 = 0$ and $\fg_1 \cdot V = V$, so $V$ is a vector space over $\F_p$; in particular it is $n!$-torsion-free and $n!$-divisible. Now Step \ref{v:oplusSymNatsl2K1:st:h} of Variation \ref{v:oplusSymNatsl2K1} requires only $n!$-torsion-freeness, so we get $h x^{n-1} = (n-1) x^{n-1}$ in $\End V$ (this is only the first step of the induction fully carried in Step \ref{v:oplusSymNatsl2K1:st:h} of Variation \ref{v:oplusSymNatsl2K1}).
As $x^{n-1} \neq 0$, we deduce $E_{n-1} \neq 0$.
Since $\oplus_{\ell \in \Z/p\Z} E_\ell$ is clearly $\fg_1$-invariant, one finds by simplicity $V = \oplus_{\ell \in \Z/p\Z} E_\ell$.

We now make the following observation: if for some $\ell \in \Z/p\Z$, $E_\ell \cap \ker x \neq 0$, then $E_\ell \leq \ker x$, and likewise with $\ker y$ instead of $\ker x$. We prove it only for $x$ as length plays no role here. Consider $W = \oplus_{i = 0}^{p-1} y^i \cdot (E_\ell \cap \ker x)$. We claim that $W$ is $x$-invariant. This is because for $a \in E_\ell \cap \ker x$ and $i \in \{0, \dots, p-1\}$ one has $x \cdot (y^0 \cdot a) = x \cdot a = 0 \in W$ when $i = 0$ and otherwise
\[x y^i \cdot a = y^i x \cdot a + i (h+i-1) y^{i-1} \cdot a \in W\]
We claim that $W$ is $y$-invariant as well. This is because for $a \in E_\ell \cap \ker x$,
\[xy\cdot (y^{p-1} \cdot a) = xy^p \cdot a = y^p x \cdot a = 0\]
whence $y \cdot (y^{p-1} \cdot a) \in E_\ell \cap \ker x \leq W$.
By assumption $W$ is non-trivial, by simplicity of $V$ one has $W = V$ and therefore $E_\ell \leq \ker x$.

Since $0 \neq \im x^{n-1} \leq \ker x \cap E_{n-1}$, one has $E_{n-1} \leq \ker x$. Now the Casimir operator $c_1 = 4yx + h^2 + 2h$ equals $n^2 - 1$ on $E_{n-1}$; by simplicity of $V$, $c_1 = n^2 - 1$ everywhere. In particular on $E_{m-1}$ one finds $4yx + m^2 - 1 = n^2 - 1 = m^2 - 1$ so $yx$ annihilates $E_{m-1}$. If $x \cdot E_{m-1} \neq 0$ then $0 \neq x \cdot E_{m-1} \leq \ker y \cap E_{m+1} = \ker y \cap E_{1-n}$, so by the above observation $y$ annihilates $E_{1-n}$ and one readily sees $V \simeq \Sym^{n-1} \Nat \fg_1$ (a very special case of our construction). If $x \cdot E_{m-1} = 0$ then one retrieves an $\bS_{\alpha, \beta}$.
\end{proofclaim}

We thus have described all simple $\fg_1$-modules of length $n$ in characteristic $\geq n+1$: they correspond to irreducible polynomials in $\F_p[X]$.
There remains one pending question: can one analyze \emph{all} $\fg_1$-modules of length $n$ in characteristic $\geq n+1$, in terms of $\bS_{\alpha, \beta}$'s? It could be conjectured so but the author wishes to dwell no longer on a subject of disputable interest.

\subsection{The Symmetric Case}\label{s:xn=yn=0}

There is something odd in assuming the characteristic of $\K_1$ to be $\geq 2n+1$ in length $n$; we bring no evidence to support the feeling that a better lower bound should be $n+1$ as it was in Variation n$^\circ$12 \cite{TV-I}.

We know from the previous subsection that lowering the characteristic can result in creating pathologies. Observe how in $\bS_{\alpha, \beta}$ the actions of $x$ and $y$ are dissymmetrical as soon as $\alpha$ or $\beta$ in non-zero. In particular $\bS_{\alpha, \beta}$ cannot be made into an $\SL_2(\K)$-module in a ``consistent'' way since $x$ and $y$ should then have the same order in $\End V$ being conjugate under the adjoint action of the Weyl group of $\SL_2(\K)$.
In short all our previous counterexamples shared the feature that the action of $y$ was quite different from that of $x$, which is ill-behaved. The minimal decency requirement on an $\sl_2(\K)$-module $V$ in order to stem from an associated $\SL_2(\K)$-module is that $x$ and $y$ should have the same order in $\End V$.

Under this extra symmetry assumption it is possible to classify $\sl_2(\F_p)$-modules of \emph{two-sided} finite length even in low characteristic.
If the author had his chance he would include the following in \cite{TV-I} rather than in the present paper.
\begin{variation}\label{v:themissingone}
Let $\K_1 = \F_3$ be the field with three elements. Let $\fg_1 = \sl_2(\F_3)$ and $V$ be a $\fg_1$-module. Suppose that $x^2 = y^2 = 0$ in $\End V$. Then $V = \Ann_V(\fg_1) \oplus \fg_1 \cdot V$, and $\fg_1\cdot V$ is a $\K_1$-vector space with $\fg_1 \cdot V \simeq \oplus_I \Nat \fg_1$ as $\K_1\fg_1$-modules.
\end{variation}
\begin{proof}
It can be computed that in $\End V$, $hx = x$ and $hy = -y$ likewise. In particular $x c_1 = x (2xy+2yx+h^2) = 2hx + x = 3 x = 0$ so $x$ annihilates $c_1 \cdot V$ and so does $y$, implying $\im c_1 \leq \Ann_V \fg_1$. Hence $c_1^2\cdot V = 0$. It can also be computed that $(h-1)h(h+1) = 0$. As $V$ may lack division by $2$ we should not jump to hasty conclusions.

Let $\overline{V} = V/\Ann_V \fg_1$. We know that $\overline{V}$ is a vector space over $\F_3$. Since the identities proved above hold in $\overline{V}$ as well, one does have $\overline{V} = E_{-1}(\overline{V}) \oplus E_0(\overline{V}) \oplus E_1(\overline{V})$.

Now $x$ annihilates $E_0(\overline{V}) \oplus E_1(\overline{V})$, as $x \cdot (E_0(\overline{V}) \oplus E_1(\overline{V})) \leq (E_{-1}(\overline{V}) \oplus E_0(\overline{V}))\cap E_1(\overline{V}) = 0$ since $\overline{V}$ has exponent $3$. Similarly, $y$ annihilates $E_{-1}(\overline{V}) \oplus E_0(\overline{V})$. It follows that $E_{0}(\overline{V}) \leq \Ann_{\overline{V}}\fg_1 = 0$ by perfectness.

Moreover $x$ is injective on $E_{-1}(\overline{V})$ since for $\overline{a_{-1}} \in E_{-1}(\overline{V}) \cap \ker x$ one has $-\overline{a_{-1}} = h \cdot \overline{a_{-1}} = (xy-yx) \cdot \overline{a_{-1}} = 0$. At this point it is clear that $\overline{V} = E_{-1}(\overline{V}) \oplus E_1(\overline{V}) \simeq \oplus_I \Nat \fg_1$.

We go back up to $V$. First observe that $\overline{E_0(V)} \leq E_0(\overline{V}) = 0$ so $E_0 (V) \leq \Ann_V \fg_1 \leq E_0(V)$.

The sum $E_{-1}(V) + E_0(V) + E_1(V)$ is direct: consider a relation $a_{-1} + a_0 + a_1$ with obvious notations. Then applying $h$ twice, one finds $-a_{-1} + a_1 = 0 = a_{-1} + a_1$. In particular $2 a_1 = 0$ and $2 (h+h)\cdot a_1 = h \cdot a_1 = a_1 = (h+h) \cdot 2 a_1 = 0$. Then $a_{-1} = 0$ and $a_0 = 0$ as well.

We claim that $V = E_{-1}(V) \oplus E_0(V) \oplus E_1(V)$. For if $v \in V$ then modulo $\Ann_V \fg_1$, $\overline{v} = \overline{v_{-1}} + \overline{v_1}$ with obvious notations, so that it suffices to lift say $\overline{v_1}$. Since $h \cdot \overline{v_1} = \overline{v_1}$ one has a priori $h\cdot v_1 = v_1 + v_0$ for some $v_0 \in E_0(V) = \Ann_V \fg_1$. Hence $h\cdot (v_1 + v_0) = v_1 + v_0$ and $v_1 = (v_1 + v_0) - v_0$ with $v_1 + v_0 \in E_1(V)$. A similar argument holds for $E_{-1}(V)$.

We finally claim that $E_{-1}(V) \oplus E_1(V)$ is $\fg_1$-invariant. If $a_1 \in E_1(V)$ then a priori $x \cdot a_1 = a_0$ with $a_0 \in \Ann_V\fg_1 = E_0(V)$. But applying $h$ one has $0 = h \cdot a_0 = h x \cdot a_1 = x \cdot a_1 = a_0$. Similarly $y \cdot a_1 = a_{-1} + a_0$ with obvious notations. Now $- a_{-1} = h y \cdot a_{-1} = - y \cdot a_{-1} = -a_{-1} - a_0$ whence $a_0 = 0$ and $y \cdot a_1 \in E_{-1}(V)$. By symmetry $E_{-1}(V) \oplus E_1(V)$ is $\fg_1$-invariant.

It is now clear that $\fg_1 \cdot V = E_{-1}(V) \oplus E_1(V) \simeq V/E_0(V) \simeq \overline{V}$ as a $\fg_1$-module. 
\end{proof}

\begin{variation}\label{v:oplusSymNatsl2K1:lowchar}
Let $n \geq 2$ be an integer and $\K_1$ be the field $\F_p$ with $n < p < 2n$.
Let $\fg_1 = \sl_2(\K_1)$ and $V$ be a $\fg_1$-module. Suppose that $x^n = y^n = 0$ in $\End V$. 

Then $V = \Ann_V(\fg_1) \oplus \fg_1 \cdot V$, and $\fg_1\cdot V$ is a $\K_1$-vector space with $\fg_1 \cdot V \simeq \oplus_{k = 1}^{n-1} \oplus_{I_k} \Sym^k \Nat \fg_1$ as $\K_1\fg_1$-modules.
\end{variation}
\begin{proof}
Induction on $n$. When $n = 2$ this is Variation \ref{v:themissingone}.
We shall adapt the proof of Variation \ref{v:oplusSymNatsl2K1}.
Write $p = n + m$ with $0 < m < n$.

One might desire to assume $\Ann_V \fg_1 = 0$ and $\fg_1\cdot V = V$. Actually if $p > n+1$ the proof given in Variation \ref{v:oplusSymNatsl2K1} remains correct as $k(k+2) \neq 0$ for $k \in \{1, \dots, n-1\}$. But when $p = n+1$ the Casimir operator $c_1$ now annihilates $\Sym^{p-2} \Nat \fg_1$ and there may be some subtleties.

\begin{step}\label{v:oplusSymNatsl2K1:lowchar:st:Fpev}
We may assume that $V$ is a $\K_1$-vector space.
\end{step}
\begin{proofclaim}
Suppose the result is known for $\K_1$-vector spaces and bear in mind that assumptions on the length of $x$ and $y$ go down to subquotients.

As $\fg_1$ annihilates $pV$ the factor $\overline{V} = V/\Ann_V\fg_1$ is a vector space. By assumption $\overline{V} = \fg_1\cdot \overline{V} \oplus \Ann_{\overline{V}} \fg_1$. Then using perfectness one has $\Ann_{\overline{V}} \fg_1 = 0$ so $\overline{V} = \fg_1 \cdot \overline{V}$.

As a consequence $V = \fg_1 \cdot V + \Ann_V \fg_1$. But $p$ annihilates the submodule $W = \fg_1 \cdot V$ which is therefore another vector space. Still by assumption $W = \fg_1 \cdot W \oplus \Ann_W \fg_1$. Then perfectness again yields $\fg_1 \cdot W = W$ so $\Ann_W \fg_1 = 0$.

In particular $W \cap \Ann_V \fg_1 = 0$ and $V = \fg_1 \cdot V \oplus \Ann_V \fg_1$ has the desired structure since $W = \fg_1 \cdot V$ does.
\end{proofclaim}

\begin{step}
In $\End V$, $(h-n+1) (h-n+2) \dots (h+n-1) = 0$.
\end{step}
\begin{proofclaim}
Since $p > n$ the proof given in Variation \ref{v:oplusSymNatsl2K1} remains correct.
\end{proofclaim}

We move to the weight space decomposition. Unfortunately the various $E_j$'s with $j \in \{1-n, \dots, n-1\}$ are no longer pairwise distinct so special attention must be paid. Observe how since $V$ is a $\K_1 = \F_p$-vector space one should actually talk about the $E_{[j]}$'s (where $[j]$ is the congruence class of $j$ modulo $p$) in order to prevent confusion. This is what we do from now on.

\begin{step}
$V = \oplus_{j \in \{0, \dots, p-1\}} E_{[j]}$.
\end{step}
\begin{proofclaim}
Bear in mind that $(h-n+1) (h-n+2) \dots (h+n-1) = 0$.
If $p = 2n-1$ the argument of Variation \ref{v:oplusSymNatsl2K1} remains correct since the polynomials $X - j$ with $j \in \{1-n, \dots, n-1\}$ are still pairwise coprime and coincide with the polynomials $X-j$ with $j \in \{0, \dots, p-1\}$. But for $p \leq 2n-3$ which we now assume it is no longer the case as some appear twice. Let us determine which with care.

As $p \leq 2n-3$ we have $n-1 \geq m+1$. We lift every congruence class modulo $p$ to its canonical representative in $\{0, \dots, p-1\}$.
\[\begin{matrix}
\mbox{class} & [1-n] & [2-n] & \dots & [-1] & [0] & [1] & \dots & [n-1]\\
\mbox{repr.} & m+1 & m+2 & \dots & p-1 & 0 & 1 & \dots & n-1
\end{matrix}\]
Let us partition $I = \{0, \dots, p-1\}$ into the set $I_1 = \{0, \dots, m\}\cup\{n, \dots, p-1\}$ of the $2m+1$ elements occurring once and the set $I_2 = \{m+1, \dots, n-1\}$ of the $n-1-m$ elements occurring twice:
\[
0, 1, \dots, m, \underbrace{m+1, \dots, n-1}_{I_2}, n, \dots, p-1
\]

Therefore the polynomial
\[P(X) = \prod_{\ell \in I_1} (X-[\ell]) \cdot \prod_{\ell \in I_2} (X-[\ell])^2\]
annihilates $h$ in $\End V$.
For $\ell \in I$ let $F_{[\ell]} = \ker (h-[\ell])^2 \geq E_{[\ell]} = \ker (h-[\ell])$. It is readily observed that $x$, resp. $y$, maps $F_{[\ell]}$ to $F_{[\ell+2]}$, resp. $F_{[\ell-2]}$.
Since all monomials powers in $P(X)$ are pairwise coprime in $\F_p[X]$ one has:
\[V = \mathop\oplus_{\ell \in I_1} E_{[\ell]} \bigoplus \mathop\oplus_{\ell \in I_2} F_{[\ell]}\]

Observe that for all $\ell \in I_1$, $F_{[\ell]} = E_{[\ell]}$. Our task is to prove it for $\ell \in I_2$ as well. So let $k \in I$ be minimal with $F_{[k]} > E_{[k]}$; $k \in I_2$ so $k \geq m+1$. We wish to take the least $i$ with $k+2i \in I_1$. Unfortunately this may fail to exist, for instance when $n=p-1$ and $k = p-2$. But there certainly is $i$ minimal with $[k+2i] \in [I_1]$. Then $i \leq \frac{n-k}{2}+1$.

Let $W = \oplus_{\ell \in I} E_{[\ell]}$ which is clearly $\fg_1$-invariant. We shall compute modulo $W$ which we denote by $\equiv$.
Let $v \in F_{[k]}$. Recall that $y \cdot v \in y \cdot F_{[k]} \leq F_{[k-2]} = E_{[k-2]} \leq W$, so $y \cdot v \equiv 0$. Moreover by construction $x^i \cdot v \in F_{[k+2i]} = E_{[k+2i]} \leq W$.
Finally by definition $(h-[k]) \cdot v \in \ker (h-[k]) = E_{[k]} \leq W$ so $h\cdot v \equiv k v$. Hence
\begin{align*}
0 & \equiv y^i x^i \cdot v\\
& \equiv \sum_{k = 0}^i (-1)^k k! \binom{i}{k} \binom{i}{k} \prod_{\ell = 0}^{k-1} (h+[\ell]) x^{i-k} y^{i-k} \cdot v\\
& \equiv (-1)^i i! \prod_{\ell = 0}^{i-1}(h+[\ell])\cdot v\\
& \equiv (-1)^i i! \prod_{\ell = 0}^{i-1} [k+\ell] v
\end{align*}
Now $k \neq 0$ since $0 \in I_1$ and $k + i - 1 \leq k + \frac{n-k}{2} \leq n < p$. Thus remains $v \equiv 0$ meaning $F_{[k]} \leq W$ and $F_{[k]} = E_{[k]}$. Therefore $V = W$.
\end{proofclaim}

\begin{localnotation}
Let $V_\bot = \im(c_1 - n^2 +1)$ and $V_\top = \ker (c_1 - n^2 + 1)$.
Let also $V_{\bot\bot} = \im(c_1 - n^2 +1)^2$ and $V_{\top\top} = \ker (c_1 - n^2 + 1)^2$.
\end{localnotation}

\begin{step}
$V_{\bot\bot}$ is a $\fg_1$-submodule isomorphic to $\oplus_{k = 1}^{n-2}\oplus_{I_k}\Sym^k \Nat \fg_1$ if $p = n+1$ and to $\Ann_{V_\bot} \oplus \oplus_{\substack{k = 1\\k \neq m-1}}^{n-2}\oplus_{I_k}\Sym^k \Nat \fg_1$ otherwise.
\end{step}
\begin{proofclaim}
As in Variation \ref{v:oplusSymNatsl2K1}, $V_\bot$ is a $\fg_1$-submodule annihilated by $x^{n-1}$, and by $y^{n-1}$ similarly. One certainly has $n-1 < p$. If $p \geq 2(n-1)+1$ then we apply Variation \ref{v:oplusSymNatsl2K1}. Otherwise $p < 2(n-1)$ and we apply induction. In any case $V_\bot = \Ann_{V_\bot} \fg_1 \oplus \fg_1 \cdot V_\bot$ and $\fg_1 \cdot V_\bot$ is isomorphic to $\oplus_{k = 1}^{n-2} \oplus_{I_k} \Sym^k \Nat \fg_1$.

The operator $c_1-n^2+1$ is no longer a bijection of $V_\bot$ as $k(k+2) = n^2 - 1$ solves into $k = n-1$ or $m-1$ in $\F_p$. Hence $c_1-n^2+1$ annihilates the component isomorphic to $\oplus_{I_{m-1}} \Sym^{m-1} \Nat \fg_1$ but acts bijectively on the other $\oplus_{I_k} \Sym^k \Nat \fg_1$'s.

As for the $\Ann_{V_\bot} \fg_1$ term, there are two possibilities. Either $p = n+1$ in which case $n^2 - 1 = 0$ and $(c_1-n^2+1)$ annihilates $\Ann_{V_\bot} \fg_1$, or $p > n+1$ in which case $n^2 - 1 \neq 0$ and $(c_1-n^2+1)$ is a bijection of $\Ann_{V_\bot} \fg_1$.

Hence if $p = n+1$ one has $V_{\bot\bot} = \oplus_{k = 1}^{n-2} \oplus_{I_k} \Sym^k \Nat \fg_1$ whereas if $p > n+1$ one has $V_{\bot\bot} = \Ann_{V_\bot} \fg_1 \oplus \oplus_{\substack{k = 1\\ k \neq m-1}}^{n-2} \oplus_{I_k} \Sym^k \Nat \fg_1$.
\end{proofclaim}

We shall simplify notations letting $\Sym^0 \Nat \fg_1$ denote the trivial $\F_p$-line so that $\Ann_{V_\bot} \fg_1$ handily rewrites into $\oplus_{I_0} \Sym^0 \Nat \fg_1$. We preferred to avoid such notation in general due to possible confusions: for instance when $V = \Z/2\Z$ as a trivial $\fg_1 = \sl_2(\F_3)$-module one has $V = \Ann_V\fg_1$ but $V$ certainly is no sum of copies of $\Sym^0 \Nat \fg_1 = \F_3$. Here we know from Step \ref{v:oplusSymNatsl2K1:lowchar:st:Fpev} that $V$ is a $\K_1 = \F_p$-vector space and confusion is no longer possible.

As a consequence $V_{\bot\bot}$ is isomorphic to $\oplus_{\substack{k = 0\\k \neq m-1}}^{n-2}\oplus_{I_k}\Sym^k \Nat \fg_1$ in either case.

\begin{step}
We may assume $V = V_{\top\top}$.
\end{step}
\begin{proofclaim}
We claim that $V = V_{\bot\bot} \oplus V_{\top\top}$. Here again $c_1-n^2+1$ is a bijection of $V_{\bot\bot}$ and its square as well whence $V_{\bot\bot} \cap V_{\top\top} = 0$ and $V = V_{\bot\bot} \oplus V_{\top\top}$.
\end{proofclaim}

From now on we suppose $V = V_{\top\top}$; in particular $(c_1 - n^2 +1)^2$ annihilates $V$.
The assumption that $y^n = 0$ in $\End V$ had played no real role up to this point.

\begin{step}
$\ker x = E_{[n-1]} \oplus E_{[m-1]}$ and $\ker y = E_{[1-n]} \oplus E_{[1-m]}$.
\end{step}
\begin{proofclaim}
We claim that $x$ is injective on $\oplus_{\substack{j \in \{0, \dots, p-1\}\\j \neq m-1, n-1}} E_{[j]}$.
Let $j \neq m-1, n-1$ and $a_j \in E_{[j]} \cap \ker x$. Then
\[0 = (c_1 - n^2 + 1)^2 \cdot a_j = (j(j+2) - n^2+1)^2 a_j\]
implies $((j+1)^2 - n^2)a_j = 0$ so by assumption on $j$ one has $a_j = 0$.

It remains to prove that $x$ does annihilate all of $E_{[n-1]} \oplus E_{[m-1]}$. First let $a_{m-1} \in E_{[m-1]}$. Do not forget that $hx^{n-1} = (n-1) x^{n-1}$ in $\End V$. So $x^{n-1}\cdot a_{m-1} \in E_{[n-3]}\cap E_{[n-1]} = 0$. But $x$ is injective on each of $E_{[1-n]}, \dots, E_{[n-5]}$ which implies $x \cdot a_{m-1} = 0$.
Hence $x$ annihilates $E_{[m-1]}$ and by symmetry $y \cdot E_{[1-m]} = 0$ as well.

Now let $a_{1-n} \in E_{[1-n]}$. Then $y^{n-1}x^{n-2} \cdot a_{1-n} \in E_{[1-n]}\cap E_{[-n-1]} = 0$. But also bearing in mind that $x$ annihilates $E_{[m-1]} = E_{[-n-1]}$:
\begin{align*}
0 = & \quad y^{n-1}x^{n-2} \cdot a_{1-n}\\
& = y \cdot \left(\sum_{k = 0}^{n-2} (-1)^k k!\ \binom{n-2}{k} \binom{n-2}{k} \left(\prod_{\ell = 0}^{k-1}(h+\ell)\right) x^{n-2-k}y^{n-2-k} \cdot a_{1-n}\right)\\
& = (-1)^{n-2} (n-2)! y \cdot (\prod_{\ell = 0}^{n-3} (h+\ell) \cdot a_{1-n})\\
& = (-1)^{n-2} (n-2)! \prod_{\ell = 0}^{n-3} (1-n+\ell) y \cdot a_{1-n}\\
& = k y \cdot a_{n-1}
\end{align*}
where $k$ is non-zero modulo $p$. Hence $y\cdot a_{1-n} = 0$. By symmetry the analogue holds of $x$.
\end{proofclaim}

We may conclude as in Variation \ref{v:oplusSymNatsl2K1}:
\[V_{\top\top} \simeq \mathop\oplus_{I_{m-1}}\Sym^{m-1}\Nat \fg_1 \bigoplus \mathop\oplus_{I_{n-1}} \Sym^{n-1}\Nat \fg_1\]
This finishes the proof.
\end{proof}

\section{Scalar Flesh}\label{S:chair}

When the irreducible $\sl_2(\K_1)$-submodules of an $\sl_2(\K)$-module $V$ are all isomorphic, $V$ bears a compatible $\K$-vector space structure: \S\ref{s:separe} contains Variation \ref{v:Symn-1Natsl2K} which is our main result.
Otherwise, and always in order to retrieve a linear geometry, one has to make some assumptions on the behaviour of $\ker x_\lambda$ and of $\im x_\lambda$. Under either assumption things work more or less in quotients of a certain composition series (\S\ref{s:suites}); should one wish to have a direct sum, one needs \emph{both} assumptions (Variation \ref{v:oplusSymNatsl2K},  \S\ref{s:separation}).

\subsection{The Separated Case}\label{s:separe}

\begin{variation}\label{v:Symn-1Natsl2K}
Let $n \geq 2$ be an integer and $\K$ be a field of characteristic $0$ or $\geq n$.
Let $\fg = \sl_2(\K)$ and $V$ be a $\fg$-module. Let $\K_1$ be the prime subfield of $\K$ and $\fg_1 = \sl_2(\K_1)$. Suppose that $V$ is a $\K_1$-vector space such that $V \simeq \oplus_I \Sym^{n-1} \Nat \fg_1$ as $\K_1\fg_1$-modules.

Then $V$ bears a compatible $\K$-vector space structure for which $V \simeq \oplus_J \Sym^{n-1} \Nat \fg$ as $\K\fg$-modules.
\end{variation}
\begin{proof}

\begin{localnotation}
For $i = 1 \dots n$ let:
\[d_i = \frac{(i-1)!\ (n-1)!}{(n-i)!} = ((i-1)!)^2 \binom{n-1}{i-1}\]
\end{localnotation}
This is an integer with prime factors $< n$. Moreover $d_{i+1} = i(n-i) d_i$.

\begin{step}\label{v:Symn-1Natsl2K:st:K1}
$V = \oplus_{i = 1}^n E_{n+1-2i}$. For all $i = 1 \dots n$, $(yx)_{|E_{n + 1 - 2 i}} = (i-1)(n + 1 - i)$, $(xy)_{|E_{n + 1 - 2 i}} = i(n - i)$, and $(x^{i-1} y^{i-1})_{|E_{n-1}} = d_i$.
\end{step}
\begin{proofclaim}
All by assumption on $V$ as a $\fg_1$-module.
\end{proofclaim}

\begin{localnotation}\label{n:scalaire}[see the linear structure in the Theme \cite{TV-I}]
Let $1 \leq i \leq n$. Set for $\lambda \in \K$ and $a_{n+1-2i} \in E_{n+1-2i}$:
\[\lambda \cdot a_{n+1-2i} = \frac{1}{n-1}\frac{1}{d_i} y^{i-1} h_\lambda x^{i-1} \cdot a_{n+1-2i} \]
\end{localnotation}

Observe that multiplication by $\lambda$ normalizes each $E_{n+1-2i}$. Extend the definition to $V = \oplus_{i = 1}^n E_{n+1-2i}$.

\begin{localremark*}
One has for all $a \in V$:
\[\lambda \cdot a = \frac{1}{n-1} \frac{1}{((n-1)!)^2} \sum_{i = 1}^n \frac{1}{d_i} y^{i-1} h_\lambda x^{n-1} y^{n-1} x^{i-1} \cdot a\]
We shall not use this.
\end{localremark*}

\begin{step}
$V$ is a $\K$-vector space.
\end{step}
\begin{proofclaim}
Let us prove that we have defined an action of $\K$. The construction is well-defined. Additivity in $\lambda$ and $a$ is obvious. So it suffices to prove multiplicativity. Let $(\lambda, \mu) \in \K^2$ and $a \in E_{n-1}$. By definition $\lambda \cdot a_{n-1} = \frac{1}{n-1} h_\lambda \cdot a_{n-1}$. So by Step \ref{v:Symn-1Natsl2K:st:K1} applied to $y_\mu \cdot a$ with $i = 2$ one has $yxy_\mu \cdot a = (n-1) y_\mu \cdot a$, whence:
\begin{align*}
(n-1)^2 \lambda \cdot (\mu \cdot a) & = h_\lambda h_\mu \cdot a\\
& = (x_\lambda y - y x_\lambda) (x y_\mu - y_\mu x) \cdot a\\
& = (x_\lambda y - y x_\lambda) x y_\mu \cdot a\\
& = x_\lambda y x y_\mu \cdot a\\
& = (n-1) x_\lambda  y_\mu \cdot a\\
& = (n-1) h_{\lambda\mu} \cdot a\\
& = (n-1)^2 (\lambda\mu) \cdot a
\end{align*}
and we obtain multiplicativity on $E_{n-1}$.

Let now $i$ be any integer in $\{1, \dots, n\}$ and $a \in E_{n+1-2i}$. Let $b = x^{i-1} \cdot a \in E_{n-1}$. Then by definition for any $\lambda \in \K$:
\[\lambda \cdot a = \frac{1}{(n-1)}\frac{1}{d_i} y^{i-1} h_\lambda \cdot b\]
so with Step \ref{v:Symn-1Natsl2K:st:K1} applied to $h_\mu \cdot b = h_\mu x^{i-1}\cdot a$:
\begin{align*}
(n-1)^2 d_i^2 \lambda \cdot (\mu \cdot a) & = y^{i-1} h_\lambda x^{i-1} y^{i-1} h_\mu x^{i-1} \cdot a\\
& = d_i y^{i-1} h_\lambda h_\mu x^{i-1} \cdot a\\
& = (n-1) d_i y^{i-1} h_{\lambda\mu} x^{i-1} \cdot a\\
& = (n-1)^2 d_i^2 (\lambda\mu) \cdot a
\end{align*}
and we obtain multiplicativity on $E_{n+1-2i}$.
\end{proofclaim}

\begin{step}
$\fg$ is linear on $V$.
\end{step}
\begin{proofclaim}
Let $\lambda \in \K$. Let us first prove linearity of $x$. It is obvious on $E_{n-1}$. So let $i \geq 2$ and $a \in E_{n+1-2i}$; one has thanks to Step \ref{v:Symn-1Natsl2K:st:K1} applied to $y^{i-2} h_\lambda x^{i-1} \cdot a \in E_{n+1-2(i-1)}$ and $i-1$:
\begin{align*}
(n-1) d_i x \cdot (\lambda \cdot a) & = x \cdot (y^{i-1} h_\lambda x^{i-1} \cdot a)\\
& = xy \cdot (y^{i-2} h_\lambda x^{i-2} \cdot (x\cdot a))\\
& = (i-1)(n+1-i) (n-1) d_{i-1} \lambda \cdot (x\cdot a)\\
& = (n-1)d_i \lambda \cdot (x\cdot a)
\end{align*}
and we obtain linearity of $x$.

Linearity of $y$ is very similar. It is obvious on $E_{n-1}$. Now for $a \in E_{n+1-2i}$ with $1 \leq i < n$ one has by Step \ref{v:Symn-1Natsl2K:st:K1} $xy \cdot a = i(n-i)a$, whence:
\begin{align*}
(n-1) d_{i+1} y \cdot (\lambda \cdot a) & = (n-1)i(n-i) d_i y\cdot (\lambda \cdot y)\\
& = y (y^{i-1} h_\lambda x^{i-1})\cdot (xy \cdot a)\\
& = y^i h_\lambda x^i \cdot (y \cdot a) \\
& = (n-1) d_{i+1} \lambda \cdot (y\cdot a)
\end{align*}
which proves linearity of $y$.

It remains to prove linearity of $h_\mu$. For $a \in E_{n-1}$ one has:
\[(n-1) h_\mu \cdot (\lambda \cdot a) = h_\mu h_\lambda \cdot a = h_\lambda \cdot (h_\mu \cdot a) = (n-1) \lambda \cdot (h_\mu \cdot a)
\]
which proves linearity of $h_\mu$ on $E_{n-1}$. Now let $i \geq 2$, $a \in E_{n+1-2i}$, and $b = x^{i-1} \cdot a \in E_{n-1}$. With Step \ref{v:Symn-1Natsl2K:st:K1} applied to $b$ one finds $y \cdot b = (n-1) x^{i-2} \cdot a$. Now remember that $x^{i-1} h_\mu = h_\mu x^{i-1} - 2 (i-1) x_\mu x^{i-2}$. Then using linearity of $x$:
\begin{align*}
(n-1) x^{i-1} h_\mu \cdot a & = (n-1) (h_\mu x^{i-1} - 2 (i-1) x_\mu x^{i-2}) \cdot a\\
& = (n-1)h_\mu \cdot b - 2 (i-1) x_\mu y \cdot b\\
& = (n-1 - 2(i-1)) h_\mu \cdot b\\
& = (n-1)(n+1-2i) \mu \cdot b\\
& = (n-1) x^{i-1} \cdot ((n+1-2i) \mu \cdot a)
\end{align*}
Since $x^{i-1}$ is injective on $E_{n+1-2i}$ by Step \ref{v:Symn-1Natsl2K:st:K1} one derives $h_\mu \cdot a = (n+1-2i) \mu \cdot a$, and this holds of any $a \in E_{n+1-2i}$. In particular by multiplicativity:
\[\lambda \cdot (h_\mu \cdot a) = \lambda \cdot ((n+1-2i) \mu \cdot a) = (n+1-2i) \mu \cdot (\lambda \cdot a) = h_\mu \cdot (\lambda \cdot a)\]
so $h_\mu$ is linear.
\end{proofclaim}
$V$ is therefore a $\K\fg$-module and its structure as such is clear. This finishes the proof.
\end{proof}

\begin{remark*}[see Variation n$^\circ$10 \cite{TV-I}]
It is now obvious that for any $\lambda \neq 0$:
$\ker x = \ker x_\lambda$, $\im x = \im x_\lambda$, $\ker y = \ker y_\lambda$, $\im y = \im y_\lambda$.
\end{remark*}

\begin{remark*}
Although our proof only requires the characteristic to be $\geq n$ it is not possible to apply the method to the modules $\bS_{\alpha, \beta}$ obtained in \S\ref{s:pathologies}. All one can get is the following which generalizes Variation n$^\circ$13 \cite{TV-I}.
\begin{quote}
{\em
Let $n \geq 2$ be an integer and $\K$ be a field of characteristic $0$ or $\geq n$.
Let $\fg = \sl_2(\K)$ and $V$ be a $\fg$-module. Let $\K_1$ be the prime subfield of $\K$ and $\fg_1 = \sl_2(\K_1)$. Suppose that $V$ is a $\K_1$-vector space such that $V \simeq \bS_{\alpha, \beta}$ as $\K_1\fg_1$-modules.

Then $V$ bears a $\K$-vector space structure such that the maps $h_\lambda$ and $x_\lambda$ are everywhere linear, but $y_\lambda$ only on $E_\ell$ for $\ell \notin \{1-n, 1-m\}$.}
\end{quote}
Preservation of the linear structure under $\alpha$ and $\beta$ depends on properties which cannot be prescribed over $\K_1$.
\end{remark*}

\subsection{Composition series}\label{s:suites}

We now prove two dual partial results.

\begin{variation}\label{v:suiteSymNatsl2K:ker}
Let $n \geq 2$ be an integer and $\K$ be a field. Let $\fg = \sl_2(\K)$ and $V$ be a $\fg$-module. If the characteristic of $\K$ is $0$ one requires $V$ to be torsion-free.
Suppose that:
\begin{itemize}
\item
either $x^n = 0$ in $\End V$ and the characteristic of $\K$ is $0$ or $\geq 2n+1$,
\item
or $x^n = y^n = 0$ in $\End V$ and the characteristic of $\K$ is $\geq n+1$.
\end{itemize}
Suppose in addition that for all $\lambda \in \K$, $\ker x \leq \ker x_\lambda$.

Then there exists a series $\Ann_V(\fg) = V_0 \leq V_1\leq \dots \leq V_{n-1} = V$ of $\fg$-submodules such that for all $k = 1, \dots, n-1$, $V_k/V_{k-1}$ bears a compatible $\K$-vector space structure for which $V_k/V_{k-1} \simeq \oplus_{I_k} \Sym^k \Nat \fg$ as $\K\fg$-modules.
\end{variation}
\begin{proof}
Induction on $n$. When $n = 2$ this is Variation n$^\circ$12 \cite{TV-I} and one even has $V = \Ann_V \fg \oplus \oplus_{I_1} \Nat \fg$.
Let $\K_1$ denote the prime subfield and $\fg_1 = \sl_2(\K_1)$.
By Variation \ref{v:oplusSymNatsl2K1} or \ref{v:oplusSymNatsl2K1:lowchar} depending on the assumptions, $V = \Ann_V(\fg_1) \oplus \fg_1 \cdot V$ where $\fg_1\cdot V \simeq \oplus_{k = 1}^{n-1} \oplus_{I_k} \Sym^k \Nat \fg_1$ as $\K_1\fg_1$-modules.

Let $V_\bot = \Ann_V \fg_1 \oplus \oplus_{k = 1}^{n-2} \oplus_{I_k} \Sym^k \Nat \fg_1$ and $V_\top =  \oplus_{I_{n-1}} \Sym^{n-1} \Nat \fg_1$.
These are $\fg_1$-submodules satisfying $V = V_\bot \oplus V_\top$.
One should be careful with the Casimir operator $c_1$. Since this operator does not commute with $\fg$ in $\End V$, $V_\bot$ and $V_\top$ as defined in Variation \ref{v:oplusSymNatsl2K1} have no reason a priori to be $\fg$-invariant. Moreover the definition of $V_\bot$ in terms of $c_1$ fails in characteristic $\leq 2n$ as seen in Variation \ref{v:oplusSymNatsl2K1:lowchar}.

Yet in the present case one sees by inspection in the $\fg_1$-module $V$:
\[V_\bot = (\oplus_{i = 1}^n E_{n-2i}) \bigoplus (\oplus_{i = 1}^n (E_{n+1-2i} \cap \ker x^{i-1}))\]

Let us now prove that $V_\bot$ is a $\fg$-submodule. It suffices to show that it is $\ft = \{h_\lambda: \lambda \in \K\}$-invariant. All $E_j$'s are $h_\lambda$-invariant. But by assumption on the kernels in $V$, $\ker x$ is as well: for if $a \in \ker x$ then $x_\lambda \cdot a = 0$ and $x h_\lambda \cdot a = (h_\lambda x - 2 x_\lambda) \cdot a = 0$. So the subgroup $V_\bot$ is $\fg$-invariant; it is a $\fg$-submodule.

\begin{localremark*}
There is no reason why $V_\top$ should be $\fg$-invariant as well.
\end{localremark*}
One sees that $x^{n-1}$ acts trivially on $V_\bot = 0$. Moreover $V_\bot$ still enjoys the property $\ker x \leq \ker x_\lambda$; induction provides the desired structure on $V_\bot$. But $V/V_\bot \simeq V_\top$ as $\fg_1$-modules so in the quotient $V/V_\bot$, $\ker x \cap \ker y^{n-1} = 0$. One then applies Variation \ref{v:Symn-1Natsl2K} to the $\fg$-module $V/V_\bot$ in order to conclude.
\end{proof}

\begin{variation}\label{v:suiteSymNatsl2K:im}
Let $n \geq 2$ be an integer and $\K$ be a field. Let $\fg = \sl_2(\K)$ and $V$ be a $\fg$-module. If the characteristic of $\K$ is $0$ one requires $V$ to be torsion-free.
Suppose that:
\begin{itemize}
\item
either $x^n = 0$ in $\End V$ and the characteristic of $\K$ is $0$ or $\geq 2n+1$,
\item
or $x^n = y^n = 0$ in $\End V$ and the characteristic of $\K$ is $\geq n+1$.
\end{itemize}
Suppose in addition that for all $\lambda \in \K$, $\im x_\lambda \leq \im x$.

Then there exists a series $0 = V_0 \leq V_1 \leq \dots \leq V_{n-1} = \fg \cdot V$ of $\fg$-submodules such that for all $k = 1, \dots, n-1$, $V_k/V_{k-1}$ bears a compatible $\K$-vector space structure for which $V_k/V_{k-1} \simeq \oplus_{I_{n-k}} \Sym^{n-k} \Nat \fg$ as $\K\fg$-modules.
\end{variation}
\begin{proof}
Induction on $n$. When $n = 2$ this is Variation n$^\circ$12 \cite{TV-I} and one even has $V = \Ann_V \fg \oplus \oplus_{I_1} \Nat \fg$.
Let $\K_1$ denote the prime subfield and $\fg_1 = \sl_2(\K_1)$.
By Variation \ref{v:oplusSymNatsl2K1} or \ref{v:oplusSymNatsl2K1:lowchar} depending on the assumptions, $V = \Ann_V(\fg_1) \oplus \fg_1 \cdot V$ where $\fg_1\cdot V \simeq \oplus_{k = 1}^{n-1} \oplus_{I_k} \Sym^k \Nat \fg_1$ as $\K_1\fg_1$-modules.

Let $V_\bot = \Ann_V \fg_1 \oplus \oplus_{k = 1}^{n-2} \oplus_{I_k} \Sym^k \Nat \fg_1$ and $V_\top =  \oplus_{I_{n-1}} \Sym^{n-1} \Nat \fg_1$.
One sees by inspection in the $\fg_1$-module $V$ that:
\[V_\top = \oplus_{i = 1}^n E_{n+1-2i} \cap \im x^{n-i}\]

Let us then prove that $V_\top$ is a $\fg$-submodule. It suffices to show that it is $\ft = \{h_\lambda: \lambda \in \K\}$-invariant. All $E_j$'s are $h_\lambda$-invariant. But by assumption on the images in $V$, $\im x$ is as well: for if $a \in \im x$ then
writing $a = x \cdot b$ one finds $h_\lambda \cdot a = x h_\lambda \cdot b + 2 x_\lambda \cdot b \in  \im x$ by assumption.
The subgroup $V_\top$ is therefore $\fg$-invariant: it is a $\fg$-submodule.

One sees that in the submodule $V_\top$, $\ker x \cap \ker y^{n-1} = 0$; Variation \ref{v:Symn-1Natsl2K} provides the desired structure on $V_\top$.
But $V/V_\top \simeq V_\bot$ as $\fg_1$-modules so in the quotient $V/V_\top$, $x^{n-1}$ acts trivially. Moreover $V/V_\top$ still enjoys the property $\im x_\lambda \leq \im x$.
One then applies induction to the $\fg$-module $V/V_\top$ in order to conclude.
\end{proof}

\subsection{Separation}\label{s:separation}

\begin{variation}\label{v:oplusSymNatsl2K}
Let $n \geq 2$ be an integer and $\K$ be a field. Let $\fg = \sl_2(\K)$ and $V$ be a $\fg$-module. If the characteristic of $\K$ is $0$ one requires $V$ to be torsion-free.
Suppose that:
\begin{itemize}
\item
either $x^n = 0$ in $\End V$ and the characteristic of $\K$ is $0$ or $\geq 2n+1$,
\item
or $x^n = y^n = 0$ in $\End V$ and the characteristic of $\K$ is $\geq n+1$.
\end{itemize}
Suppose in addition that for all $\lambda \in \K$, $\ker x \leq \ker x_\lambda$ and $\im x_\lambda \leq \im x$.

Then $V = \Ann_V(\fg) \oplus \fg \cdot V$, and $\fg\cdot V$ bears a compatible $\K$-vector space structure for which $\fg \cdot V \simeq \oplus_{k = 1}^{n-1} \oplus_{I_k} \Sym^{n-1} \Nat \fg$ as $\K\fg$-modules.
\end{variation}
\begin{proof}
Induction on $n$. When $n = 2$ this is Variation n$^\circ$12 \cite{TV-I}.
As in Variations \ref{v:suiteSymNatsl2K:ker} and \ref{v:suiteSymNatsl2K:im}, $V_\bot$ and $V_\top$ are $\fg$-invariant.
But the property $\ker x \leq \ker x_\lambda$ clearly goes to submodules, and the property $\im x_\lambda \leq \im x$ clearly goes to quotients. Hence $V_\bot \simeq V/V_\top$ (here as $\fg$-modules) allows to use induction.
\end{proof}

\section{Lesson: coherence degrees}\label{S:coherence}

\begin{notation*}
Let $V$ be a $\fg$-module. Let:
\begin{itemize}
\item
$\kappa(V)$\inmargin{$\kappa(V)$} be the least integer $n$, if there is one, such that for all $(\lambda_1, \dots, \lambda_n) \in \K^n$, $\ker x^n \leq \ker x_{\lambda_1}\dots x_{\lambda_n}$;
\item
$\iota(V)$\inmargin{$\iota(V)$} be the least integer $n$, if there is one, such that for all $(\lambda_1, \dots, \lambda_n) \in \K^n$, $\im x_{\lambda_1}\dots x_{\lambda_n} \leq \im x^n$.
\end{itemize}
\end{notation*}


The parameters $\kappa(V)$ and $\iota(V)$ may play some role in the rest of the present series of articles. A convenient name would be the \emph{ascending} (resp., \emph{descending}) \emph{coherence degrees} of the action. Be careful that they are \emph{not} the least $n$ such that the kernels (resp. images) of $x_{\lambda_1} \dots x_{\lambda_n}$ do not depend on $(\lambda_1, \dots, \lambda_n)$. They are the least $n$ such that one always has an inclusion.

We have in Variation \ref{v:oplusSymNatsl2K} been using an obvious fact.
\begin{observation*}
Let $V$ be a $\fg$-module and $W \leq V$ be a $\fg$-submodule. Then
$\kappa(W) \leq \kappa(V)$ and $\iota(V/W) \leq \iota(V)$.
\end{observation*}

Remember that $\lambda(V)$ stands for the length of $V$ as a $\fu$-module. One knows from  Variation \ref{v:xn=0l=0} that $\lambda(V)$ is the length of $V$ as an $x$-module, at least provided the characteristic is not too low.

\begin{variation}\label{v:xn=0n-1coherente}
Let $n \geq 2$ be an integer and $\K$ be a field of characteristic $0$ or $\geq n+1$.
Let $V$ be a $\sl_2(\K)$-module of $\fu$-length at most $n$.
Then for all $\lambda_1, \dots, \lambda_{n-1} \in \K$ one has $\ker x^{n-1} \leq \ker (x_{\lambda_1} \dots x_{\lambda_{n-1}})$ and $\im (x_{\lambda_1} \dots x_{\lambda_{n-1}}) \leq \im x^{n-1}$.
\end{variation}
In our notations this writes $\kappa(V) \leq \lambda(V) - 1$ and $\iota(V) \leq \lambda(V) - 1$.
\begin{proof}
Let us first deal with the kernels. We shall need the following identity of the enveloping ring (remember that the terms in the hats do not appear):
\begin{align}
x_{\lambda_1} \dots x_{\lambda_i} y_\mu & = \quad y_\mu x_{\lambda_1} \dots x_{\lambda_i} + \sum_{j } h_{\mu \cdot \lambda_j} x_{\lambda_1} \dots \widehat{x_{\lambda_j}} \dots x_{\lambda_i}\nonumber\\
& \quad - \sum_{j \neq k} x_{\mu\cdot \lambda_j \cdot \lambda_k}x_{\lambda_1} \dots \widehat{x_{\lambda_j}} \dots \widehat{x_{\lambda_k}} \dots x_{\lambda_i}\label{x_ulambday_mu}
\end{align}

We prove by induction on $i = 0 \dots n-1$ that for all $(\lambda_1, \dots, \lambda_i) \in \K^i$, $x^{n-1-i} x_{\lambda_1} \dots x_{\lambda_i}$ annihilates $\ker x^{n-1}$.
When $i = 0$ this is obvious. Let us suppose that the property holds of $i$ and prove it of $i+1\leq n-1$. Let $(\lambda_1, \dots, \lambda_i, \mu)$ be an $(i+1)$-tuple of $\K$ and set $\pi = x^{n-1 - (i+1)} x_{\lambda_1} \dots x_{\lambda_i} x_\mu$.
Recall that $2 x_\mu = 2 x y_\mu x - y_\mu x^2 - x^2 y_\mu$.
By assumption on the length all products of the form $x^{n-i} x_{\lambda_1} \dots x_{\lambda_i}$ are zero, whence in $\End V$:
\begin{align*}
2 \pi & = 2 x^{n-1 - (i+1)} x_{\lambda_1} \dots x_{\lambda_i} x_\mu\\
& = x^{n-2-i} x_{\lambda_1} \dots x_{\lambda_i} (2 x y_\mu x - y_\mu x^2 - x^2 y_\mu)\\
& = 2 x^{n-1-i} x_{\lambda_1} \dots x_{\lambda_i} y_\mu x - x^{n-2-i} x_{\lambda_1} \dots x_{\lambda_i} y_\mu x^2
\end{align*}

It remains to move the $y_\mu$'s to the left using equation (\ref{x_ulambday_mu}) applied to the tuples $(x, \dots, x, x_{\lambda_1}, \dots, x_{\lambda_i})$. Let us do it mentally.
Terms with a $y_\mu$ on the left will end in $x^{n-i} x_{\lambda_1} \dots x_{\lambda_i}$: by assumption they are zero in $\End V$. Terms with a $h_\nu$ on the left end either in $x^{n-1-i} x_{\lambda_1} \dots x_{\lambda_i}$ or in $x^{n-i} x_{\lambda_1} \dots \widehat{x_{\lambda_j}} \dots x_{\lambda_i}$ for some $j$: by induction they annihilate $\ker x^{n-1}$. It thus only remains to consider the pure products of $x$ and the various $x_\nu$'s. There are three cases:
\begin{itemize}
\item
the $j\th$ and $k\th$ (omitted) terms were among the $x_\nu$'s: the product is of the form $x_{\mu\cdot \nu_1 \cdot \nu_2} x^{n-i} x_{\lambda_1} \dots \widehat{x_{\nu_1}} \dots \widehat{x_{\nu_2}} \dots x_{\lambda_i}$; by induction it annihilates $\ker x^{n-1}$.
\item
the $j\th$ (omitted) term was among the $x$'s and the $k\th$ among the $x_\nu$'s (or vice-versa): the product is of the form $x_{\mu\cdot \nu} x^{n-1-i} x_{\lambda_1} \dots \widehat{x_\nu} \dots x_{\lambda_i}$; it annihilates $\ker x^{n-1}$.
\item
the $j\th$ and $k\th$ (omitted) terms were among the $x$'s: the product is of the form $x_\mu x^{n-i-2} x_{\lambda_1} \dots x_{\lambda_i} = \pi$.
\end{itemize}
The latter case is of interest. Paying attention to the signs and coefficients it appears exactly $-4\binom{n-1-i}{2} + 2\binom{n-2-i}{2} = - (n-2-i)(n+1-i)$ times, whence:
\[2 \pi = - (n-2-i)(n+1-i) \pi + z\]
where $z$ annihilates $\ker x^{n-1}$, that is $(n-i-1)(n-i) \pi$ annihilates $\ker x^{n-1}$. Now $i \leq n-2$ so if we had started with $\frac{\mu}{(n-i-1)(n-i)}$ we would have found that $\pi$ annihilates $\ker x^{n-1}$. This completes the induction; with $i = n-1$ one obtains the desired conclusion.

As far as the images are concerned we proceed similarly using the dual formula:
\begin{align}
y_\mu x_{\lambda_1} \dots x_{\lambda_i} & = \quad x_{\lambda_1} \dots x_{\lambda_i} y_\mu - \sum_{j } x_{\lambda_1} \dots \widehat{x_{\lambda_j}} \dots x_{\lambda_i} h_{\mu \cdot \lambda_j} \nonumber\\
& \quad - \sum_{j \neq k} x_{\lambda_1} \dots \widehat{x_{\lambda_j}} \dots \widehat{x_{\lambda_k}} \dots x_{\lambda_i} x_{\mu\cdot \lambda_j \cdot \lambda_k}\label{y_mux_ulambda}
\end{align}
proving by induction on $i = 0\dots n-1$ that for all $(\lambda_1, \dots, \lambda_i)\in \K^i$ one has $\im (x^{n-1-i} x_{\lambda_1} \dots x_{\lambda_i}) \leq \im x^{n-1}$. When rewriting $\pi$ use instead:
\[2 \pi = (2 x y_\mu x - y_\mu x^2 - x^2 y_\mu) x^{n-2-i} x_{\lambda_1} \dots x_{\lambda_i}\]
and move the $y_\mu$'s to the right using formula (\ref{y_mux_ulambda}).
\end{proof}

\begin{remark*}
The author does not feel comfortable with writing a double proof; there is an obvious redundancy.
Some inner duality must be present but I cannot see which.
\end{remark*}

\begin{remark*}\
\begin{itemize}
\item
Equalities may not hold in Variation \ref{v:xn=0n-1coherente}: remember that in $\Nat \sl_2 (\C) \otimes{}^\varphi \Nat \sl_2 (\C)$ ($\varphi$ stands for complex conjugation) one has $x^3 = 0$, $x_i x_1 = 0$, and $x^2 \neq 0$.
\item
The value $n-1$ is optimal. Take distinct field automorphisms $\varphi_1, \dots, \varphi_n$ and set $V = ({}^{\varphi_1} \Nat \sl_2) \otimes \dots \otimes ({}^{\varphi_n} \Nat \sl_2)$. This is an irreducible representation. Its length is $n+1$; in particular $\ker x^n \leq \ker x_{\lambda_1}\dots x_{\lambda_n}$ for all $(\lambda, \dots, \lambda_n) \in \K^n$, but this fails at stage $n-1$.

Let indeed $\lambda \in \K$ be such that $\varphi_1 (\lambda) \neq \varphi_n(\lambda)$. The standard basis $(e_1, e_2)$ of $\Nat \sl_2$ being fixed, $e_{i_1, \dots, i_n}$ will denote the pure tensor $e_{i_1} \otimes \dots \otimes e_{i_n}$. Consider $a = e_{2, \dots, 2, 1} - e_{1, 2, \dots, 2}$; one sees that $x^{n-1} \cdot a = 0$ but $x_\lambda x^{n-2} \cdot a = (n-2)! (\varphi_1(\lambda) - \varphi_n(\lambda)) e_{1, \dots, 1} \neq 0$.
\end{itemize}
\end{remark*}

One might expect $\kappa(V)$ and $\iota(V)$ to provide an indication of the number of tensor factors; but one would first need to conjecture that every simple $\fg$-module of finite length is a tensor product of copies, twisted by field automorphisms, of a same representation  of $\fg$ as a Lie \emph{algebra}. The author does not wish to do so even under model-theoretic assumptions. Anyway we have until now been dealing mostly with actions of coherence degree $1$, in a sense or the other.

It is not a priori clear that $\kappa(V)$ and $\iota(V)$ need in general be equal and the question deserves to be asked, at least for an action of finite length. Note that one could define the same numbers for the action of $y$; perhaps one should not expect a relation to the coherence degrees for $x$ even in the finite length case.

Finally, an alternative indicator could be the nilpotence height of the Casimir operator, that is the least $n$ such that $[\fg, \dots, [\fg, c_1]]$ acts trivially on $V$. Our results would have been more naive under the assumption that $c_1$ commutes with the action of $\fg$ since instead of Variations \ref{v:suiteSymNatsl2K:ker}, \ref{v:suiteSymNatsl2K:im}, and \ref{v:oplusSymNatsl2K} it would have sufficed to adapt the rather standard techniques of Variation \ref{v:oplusSymNatsl2K1}. Besides we found no relation between the nilpotence height of the Casimir operator and the coherence degrees.

One easily imagines how to define $\lambda, \kappa, \iota$ for an action of $G = \SL_2(\K)$.
\medskip
\hrule
\medskip
\noindent
Future variations will explore the symmetric powers of $\Nat G$.

\bibliographystyle{plain}
\bibliography{Variationen}

\end{document}